\documentclass[12pt,reqno]{amsart}
\usepackage{latexsym}
\usepackage{amsfonts}
\usepackage{amsmath}
\usepackage{amsthm}
\usepackage{amssymb}
\usepackage{amscd}
\usepackage[dvips]{graphicx}
\usepackage{eepic}
\def\Pa{\overline{P}}
\def\Pb{\widehat{P}}

\def\inv{\mathop {\rm inv}}
\def\maj{\mathop {\rm maj}}

\newcommand\rob{\mathop{\rm rob}}

\newcommand\ros{\mathop{\rm ros}}
\newcommand\rcb{\mathop{\rm rcb}}
\newcommand\lcb{\mathop{\rm lcb}}
\newcommand\lob{\mathop{\rm lob}}
\newcommand\rcs{\mathop{\rm rcs}}
\newcommand\lcs{\mathop{\rm lcs}}

\newcommand\rsb{\mathop{\rm rsb}}

\newcommand\los{\mathop{\rm los}}

\newcommand\mak{\mathop{\rm mak}}
\newcommand\lmak{\mathop{\rm lmak}}
\newcommand\bmaj{\mathop{ \rm bmaj}}

\newcommand\cbmaj{\mathop{ \rm cbMaj}}

\newcommand\stat{\mathop{ \rm stat}}

\newcommand\op{\mathop{ \rm op}}
\newcommand\cl{\mathop{ \rm cl}}

\def\O{{\mathcal O}}
\def\N{\mathbb{N}}
\def\S{{\mathcal S}}
\def\P{{\mathcal P}}
\def\OP{{\mathcal O}{\mathcal P}}

\newtheorem{thm}{Theorem}[section]
\newtheorem{prop}[thm]{Proposition}

\newtheorem{cor}[thm]{Corollary}

\newtheorem{lem}[thm]{Lemma}
\newtheorem{conj}[thm]{Conjecture}

\newtheorem{rmk}[thm]{Remark}

\def\pf{\noindent {\it Proof.} }

\makeatletter
\newfont{\footsc}{cmcsc10 at 8truept}
\newfont{\footbf}{cmbx10 at 8truept}
\newfont{\footrm}{cmr10 at 10truept}
\makeatother 

\newtheorem{corollary}[thm]{Corollary}
\newtheorem{defn}[thm]{Definition}
\newtheorem{lemma}[thm]{Lemma}

\numberwithin{equation}{section}
\newenvironment{demo}[1]{%
  \trivlist
  \item[\hskip\labelsep
        {\bf #1.}]
}{%
\hfill\qedsymbol
  \endtrivlist
}
\newcommand\euler[2]{\left\langle{#1}\atop {#2}\right\rangle}
\def\open{\rm{open}}
\def\clos{\rm{clos}}
\def\block{\rm{block}}
\def\mak{\rm{mak}}
\def\lmak{\rm{lmak}}
\def\los{\rm{los}}
\def\lcs{\rm{lcs}}
\def\ros{\rm{ros}}
\def\rcs{\rm{rcs}}
\def\lsb{\rm{lsb}}
\def\rsb{\rm{rsb}}
\def\lob{\rm{lob}}
\def\lcb{\rm{lcb}}
\def\rob{\rm{rob}}
\def\rcb{\rm{rcb}}
\def\bInv{\rm{bInv}}
\def\exc{\rm{cinv}}
\def\bExc{\rm{bExc}}
\def\cbInv{\rm{cbInv}}
\def\bMaj{\rm{bMaj}}

\def\cinvLSB{\rm{cinvLSB}}
\def\cmajLSB{\rm{cmajLSB}}
\def\OP{\mathcal{OP}}

\def\O{\mathcal{O}}
\def\C{\mathcal{C}}
\def\S{\mathcal{S}}
\def\T{\mathcal{T}}
\def \form{\mathcal{F}}

\newcommand\perm{\operatorname{perm}}

\newcommand\sgn{\operatorname{sgn}}

\newcommand\rdots{\mathinner{\mkern1mu\raise0pt\vbox{\kern7pt\hbox{.}}
     \mkern2mu\raise4pt\hbox{.}\mkern2mu\raise8pt\hbox{.}\mkern1mu}}

\def\qbinom#1#2{{\left[{{#1}\atop{#2}}\right]_{q}}}
\def\UB{{\overline N}}
\def\LB{{\widehat N}}
\def\UUB{{\check N}}
\def\DB{{N'}}
\def\UM{{\overline M}}
\def\LM{{\widehat M}}
\def\UUM{{\check M}}

\def\N#1{{\widehat{#1}}}
\def\K#1{{\operatorname{K}_{#1}}}
\def\vecX{{X}}
\def\vecU{{X'}}
\def\vecV{{X''}}
\def\UB{{\overline N}}
\def\LB{{\widehat N}}
\def\UUB{{\check N}}
\def\DB{{\dot{N}}}
\def\UM{{\overline M}}
\def\LM{{\widehat M}}
\def\UUM{{\check M}}

\def\N#1{{\widehat{#1}}}
\def\K#1{{\operatorname{K}_{#1}}}
\def\vecX{{X}}
\def\vecU{{\dot{X}}}
\def\vecV{{\ddot{X}}}
\def\Col#1{{\operatorname{Col}_{#1}}}
\def\entryB{{\dot{b}}}
\def\seqpoly{F}
\def\seqbinom#1#2{{\left[{{#1}\atop{#2}}\right]_{\seqpoly}}}

\begin{document}
\title{Statistics on Ordered
Partitions of Sets and $q$-Stirling Numbers}~\thanks{Version of 06/06/2006}
\author{Masao Ishikawa}
\address{Faculty of Education, Tottori University\\
\small Koyama, Tottori, Japan} \email{ishikawa@fed.tottori-u.ac.jp}

\author{Anisse Kasraoui}
\address{Institut Camille Jordan,
        Universit\'e Claude Bernard (Lyon I)\\
        F-69622, Villeurbanne Cedex, France}
\email{anisse@math.univ-lyon1.fr}

\author{Jiang Zeng}
\address{Institut Camille Jordan,
        Universit\'e Claude Bernard (Lyon I)\\
        F-69622, Villeurbanne Cedex, France}
\email{zeng@math.univ-lyon1.fr}

\begin{abstract}
An ordered partition of $[n]:=\{1,2,\ldots, n\}$ is a sequence of its
disjoint subsets whose union is $[n]$. The number of ordered partitions of $[n]$
with $k$ blocks is
$k!S(n,k)$, where $S(n,k)$ is the Stirling number of second kind.
 In this paper we prove some refinements of this formula by showing that
  the generating function of some statistics on the set of ordered partitions
  of $[n]$ with $k$ blocks is a natural $q$-analogue of $k!S(n,k)$.
In particular,  we prove several conjectures of  Steingr\'{\i}msson.
To this end, we construct a mapping  from ordered
partitions to walks in some digraphs and
 then, thanks to transfer-matrix method, we determine the
  corresponding generating functions by determinantal computations.

\end{abstract}
\maketitle

{\small \tableofcontents }

 \vskip 0.5cm \noindent{\bf MR Subject
Classifications}: Primary 05A18; Secondary 05A15, 05A30.

\section{Introduction}
An \emph{(unordered) set partition} of $[n]=\{1,2,\ldots, n\}$ is
a collection of its disjoint subsets, called $\emph{blocks}$, whose union
is $[n]$.  By convention, the standard notation of a partition of
$[n]$ is $\pi_0=B_1-B_2-\cdots -B_k$, where the blocks $B_i$
 are arranged  in increasing order of
their minimal elements and in each block $B_i$ the elements are
arranged in increasing order. Let $|\pi|=n$ if $\pi$ is a
partition of $[n]$. Let $\P_n^k$ be the set of partitions of $[n]$.

An \emph{ordered partition} $\pi$ of $[n]$ with $k$ blocks
is a rearrangement of blocks of a partition in $\P_n^k$.
Namely
$\pi=B_{\sigma(1)}-B_{\sigma(2)}-\cdots-B_{\sigma(k)}$,
where $\sigma$ is a permutation of $[k]$. We will say that $\sigma$ is the
 permutation induced by $\pi$ and set $\sigma=\perm(\pi)$.

Let $\OP_{n}^k$ be the set of ordered partitions of $[n]$ into $k$
blocks, $\OP_n = \bigcup_{k\geq 1}\OP_{n}^k$ be the set of all
ordered partitions of $[n]$, and $\OP^k= \bigcup_{n\geq 1}\OP_{n}^k$
be the set of all ordered partitions into $k$ blocks.
Clearly we
have $|\OP_n|=\sum_{k=0}^nk!S(n,k)$, where $S(n,k)$ is the Stirling  number of
second kind and it is not hard to derive
 the following exponential generating function~:
\[
\sum_{n\geq 0}\,|\OP_{n}|\frac{z^n}{n!} =\frac{1}{2-e^z}=
1+z+3\,\frac{z^2}{3!}+13\,\frac{z^3}{3!}+75\, \frac{z^4}{4!}+\cdots.
\]

Define the $p,q$-integer $[n]_{p,q}=\frac{p^n-q^n}{p-q}$, the
$p,q$-factorial $[n]_{p,q}!=[1]_{p,q}[2]_{p,q}\cdots[n]_{p,q}$ and
the $p,q$-binomial coefficients
$$
{n\brack k}_{p,q}=\frac{[n]_{p,q}!}{[k]_{p,q}![n-k]_{p,q}!} \qquad
n\geq k\geq 0.
$$
If $p=1$, we shall write $[n]_q$,  $[n]_q!$ and ${n\brack k}_q$
 for $[n]_{1,q}$,  $[n]_{1,q}!$ and ${n\brack k}_{1,q}$
respectively.

The following  $q$-analogues of \emph{Eulerian numbers}
and \emph{Stirling numbers of the second kind} were first introduced by
Carlitz~\cite{Ca1,Ca2}.

 The $q$-Eulerian numbers $\euler{n}{k}_q$ ($n\geq k\geq 0$) are
defined by
$$
\euler{n}{k}_q=q^k[n-k]_q\euler{n-1}{k-1}_q
+[k+1]_q\euler{n-1}{k}_q.
$$
The first values of the $q$-Eulerian numbers $\euler{n}{k}_q$
($n\geq k\geq 0$) read
$$
\begin{tabular}{c| ccccc}
$n\setminus k$&0&1&2&3&\\
\hline
1&1&&&&\\
2& 1&q&&&\\
3&1&$2q+2q^2$&$q^3$&&\\
4&1&$3q+5q^2+3q^3$&$3q^3+5q^4+3q^5$&$q^3$.
\end{tabular}
$$

Let $\sigma=\sigma(1)\sigma(2)\ldots \sigma(n)$ be a permutation of
$[n]$, the integer $i\in [n-1]$ is called a \emph{descent}
 of $\sigma$ if $\sigma(i)>\sigma(i+1)$. The \emph{major index} of $\sigma$,
 noted $\maj\sigma$,
 is the sum of its descents, i.e.,
$\maj\sigma=\sum_i i$, where the summation is over all descents $i$ of $\sigma$.
Then Carlitz~\cite{Ca2} gave the following combinatorial interpretation of $q$-Eulerian numbers:
$$
\euler{n}{k}_q=\sum_{\sigma}q^{\maj\, \sigma},
$$
where the summation is over all permutations of $[n]$  with $k$
descents.

The $q$-Stirling numbers  $S_q(n,k)$ of the second kind  are defined
by:
\begin{eqnarray}\label{eq:stirling}
 S_q(n,\, k)=q^{k-1} S_q(n-1,\,k-1)+[k]_qS_q(n-1,\, k)\qquad (n\geq k\geq 0),
\end{eqnarray}
 where $S_q(n,k)=\delta_{n\,k}$ if $n=0$ or $k=0$.
The first values of the $q$-Stirling numbers $S_q(n,\, k)$  read
$$
\begin{tabular}{c|ccccc}
$n\setminus k$&1&2&3&4&\\
\hline
1&1&&&&\\
2&1&q&&&\\
3&1&$1+q+q^2$&$q^3$&&\\
4&1&$1+3q+2q^2+q^3$&$q^2+2q^3+2q^4+q^5$&$q^6$.
\end{tabular}
$$

There has been a considerable amount of recent interest in
properties and combinatorial interpretations of the $q$-Eulerian numbers and $q$-Stirling
numbers and related numbers (see e.g.
\cite{Ca1,Ca2,CSZ,KsZe,Mi,Re,Sa,SS1,SS2,Stein,Wa,WW,Wh}).

The following identity  was derived in \cite{ZZ}:
\begin{equation}\label{eq:zezh}
[k]_q!\,S_q(n,k)=\sum_{m=1}^{k}q^{k(k-m)}\,{n-m\brack
n-k}_q\,\euler{n}{m-1}_q.
\end{equation}
 In the aim to
give a combinatorial proof of \eqref{eq:zezh},
Steingr\'{\i}msson~\cite{Stein} introduced the following
\begin{defn}
 A statistic $Stat$ on $\OP_n^k$ is called {\sl Euler-Mahonian} if its generating
function  is equal to $[k]_q!S_q(n,k)$, i.e.,
$$
\sum_{\pi \in \OP_n^k}q^{Stat\,\pi}=[k]_q!\,S_q(n,k).
$$
\end{defn}

Steingr\'{\i}msson~\cite{Stein} has found a few of {\sl Euler-Mahonian}
 statistics and
conjectured more  such statistics on ordered partitions. From
a different point view, Wachs~\cite{Mi} has also obtained some
 {\sl Euler-Mahonian} statistics on
ordered partitions.
Although Zeng~\cite{Ze} has showed that
much more such statistics can be derived from some
classical bijections between ordered partitions and weighted Motzkin paths,
it is not clear how to encode the conjectured statistics of Steingr\'{\i}msson
by the statistics obtained by this method.

It is the purpose of this paper to propose a new approach to attack such kind
of problem. We shall
construct a bijection $\psi$ between ordered partitions and
some walks in some
digraphs (see section 3).  This bijection
keeps track of several statistics of Steingr\'{\i}msson. Then, by transfer-matrix
method, we evaluate the generating functions of these statistics on ordered
partitions and prove that they are indeed Euler-Mahonian.


\section{Definitions and main results}

\subsection{Definitions}
Let $\pi=B_1-B_2-\cdots-B_k$ be a partition in $\OP_n^k$.
  The \emph {opener} of a block in $\pi$ is its least element
  and the \emph {closer} is its greatest element.
  The sets of openers and
closers of $\pi$ are denoted by $\open (\pi)$ and $\clos (\pi)$, respectively.
 We define a \emph{partial
order} on blocks $B_i's$ as follows~: $B_i>B_j$ if all the letters
of $B_i$ are greater than those of $B_j$; in other words, if the
opener of $B_i$ is greater than the closer of $B_j$.
 We say that $i$ is a
\emph{block descent} in $\pi$ if $B_i>B_{i+1}$. The \emph{block
major index} of $\pi$, denoted $\bmaj(\pi)$, is the sum of the block
descents in $\pi$.
  A \emph{block excedance} (resp. \emph{block
inversion}) in $\pi$ is a pair $(i,j)$ such that $i<j$ and
$B_{i}<B_{j}$ (resp. $B_{i}>B_{j}$). We denote by $\bExc \pi$ (resp.
$\bInv \pi$) the number of block excedances (resp. block inversions)
in $\pi$.
 Let \block(i) be the index
of the block (counting from the left) containing $i$, namely the
integer $j$ such that $i\in B_j$.

Following Steingr\'{\i}msson~\cite{Stein},
for $1\leq i\leq k$
we define ten coordinate statistics on $\pi \in \OP_{n}^{k}$~:
\begin{eqnarray*}
\ros_i(\pi) &=& \# \{j \in \open(\pi)\, |\,i>j, \,\block(j)>\block(i)\}, \\
\rob_i(\pi) &=& \# \{j \in \open(\pi)\, |\,i<j, \,\block(j)>\block(i)\}, \\
\rcs_i(\pi) &=& \# \{j \in \clos(\pi)\, |\,i>j, \,\block(j)>\block(i)\},\\
\rcb_i(\pi) &=& \# \{j \in \clos(\pi)\, |\,i<j, \,\block(j)>\block(i)\},\\
\los_i(\pi) &=& \# \{j \in \open(\pi)\, |\,i>j, \,\block(j)<\block(i)\},\\
\lob_i(\pi) &=& \# \{j \in \open(\pi)\, |\,i<j, \,\block(j)<\block(i)\}, \\
\lcs_i(\pi) &=& \# \{j \in \clos(\pi)\, |\,i>j, \,\block(j)<\block(i)\}, \\
\lcb_i(\pi) &=& \# \{j \in \clos(\pi)\, |\,i<j,\,\block(j)<\block(i)
\},
\end{eqnarray*}
and let $\rsb_i(\pi)$ (resp. $\lsb_i(\pi)$) be the number of
blocks \textrm{B} in $\pi$ to the
right (resp. left) of the block
containing $i$ such that the opener of \textrm{B} is smaller than
$i$ and the closer of \textrm{B} is greater than $i$.
Then define \ros, \rob, \rcs, \rcb, \lob, \los,  \lcs, \lcb, \lsb \; and\;  \rsb \; as the sum
of their coordinate statistics, e.g.
$$
\ros=\sum_{i}\ros_{i}.
$$
For any set of nonnegative integers $A$ and a composed statistic
$Stat$ on ordered partitions, we define $stat(A)$ as the sum of the
coordinate statistics in A, i.e.,
  $$
  \stat(A):=\sum_{i\in A}stat_i.
  $$
Now, for any mapping $f$ from $\OP_n^k$ to the set of subsets of $[n]$,
we define
$stat(f)$ by $stat(f)(\pi):=stat(f(\pi))$.

For a permutation $\sigma$ of $[n]$, the
pair $(i,j)$ is an \emph{inversion} if  $1\leq
i<j\leq n$ and $\sigma(i)>\sigma(j)$.
Let $\inv\sigma$ be the number of inversions in $\sigma$ and
$$\exc\sigma={n\choose 2}-\inv\sigma.
$$
By convention, for a partition $\pi$,
 we put $\inv\,\pi=\inv(\perm(\pi))$ and $\exc\,\pi=\exc(\perm(\pi))$.
Note that  $\P^k=\{\pi\in\OP^k\:|\:\inv\,\pi=0\}$ and  $\bInv\,\pi=0$ for each $\pi\in \P^k$.

Given an ordered partition $\pi$, let $\pi^r$ be the ordered
partition obtained from $\pi$ by reversing the order of the blocks.
This turns a left (resp. right) opener  into a right (resp. left)
opener, and likewise for the closers.
 Moreover, let $\pi^c$ be the ordered partition obtained by
complementing each of the letters in $\pi$, that is , by replacing
the letter $i$ by $n+1-i$. Then, it is easy to see that
$\rob(\pi^c)=\rcs(\pi)$ and $\ros(\pi^c)=\rcb(\pi)$, and likewise
for the left and closer statistics.
 Thus the eight statistics obtained by independently varying
left/right, opener/closer and smaller/bigger
 fall into only two categories when it comes to their distribution
on ordered partitions. One of these categories consists of $\rob,
\lob, \rcs$ and $\lcs$, and the other contains $\ros, \los,
\rcb$ and $\lcb$. Note that these results are completely false on
the unordered set partitions.

 For instance, we give the values of the coordinate statistics
computed
 on the partition $\pi=  6\;8  -  5  -  1\;4\;7  -  3\;9
- 2$ :
$$
\begin{array}{cccccccccc}
\pi=&  6\;8 & - & 5 & - & 1\;4\;7 & - & 3\;9 & - & 2\\
&&&&&&&&&\\
   \los_i: & 0\;0 & - & 0 & - & 0\;0\;2 & - & 1\;3 & - & 1 \\
   \ros_i: & 4\;4 & - & 3 & - & 0\;2\;2 & - & 1\;1 & - & 0 \\
   \lob_i: & 0\;0 & - & 1 & - & 2\;2\;0 & - & 2\;0 & - & 3 \\
   \rob_i: & 0\;0 & - & 0 & - & 2\;0\;0 & - & 0\;0 & - & 0 \\
   \lcs_i: & 0\;0 & - & 0 & - & 0\;0\;1 & - & 0\;3 & - & 0 \\
   \rcs_i: & 2\;3 & - & 1 & - & 0\;1\;1 & - & 1\;1 & - & 0 \\
   \lcb_i: & 0\;0 & - & 1 & - & 2\;2\;1 & - & 3\;0 & - & 4 \\
   \rcb_i: & 2\;1 & - & 2 & - & 2\;1\;1 & - & 0\;0 & - & 0 \\
   \lsb_i: & 0\;0 & - & 0 & - & 0\;0\;1 & - & 1\;0 & - & 1 \\
   \rsb_i: & 2\;1 & - & 2 & - & 0\;1\;1 & - & 0\;0 & - & 0 \\
\end{array}
$$

 Note that there are four block inversions: {\small $\{6,8\}>\{5\},\; \{6,8\}>\{2\},\;
\{5\}>\{2\}$} and $\{3,9\}>\{2\}$, and two block descents at $i=1$
and $4$; thus $\bInv\,\pi=4$ and $\bmaj\,\pi=1+4=5$. Note also that
$\bExc\,\pi=0$.
Moreover, $perm(\pi)=54132$ and thus $\inv(\pi)=8$ and $\exc(\pi)={5\choose2}-8=2$. \\

 Inspired by a statistic $\mak$ due to  Foata $\&$
 Zeilberger~\cite{FoZe}
on the  permutations, Steingr\'{\i}msson introduced its analogous on
$\OP_n^k$ as follows:
\begin{eqnarray*}
\mak&=& \ros +\lcs,\\
\lmak{}&=& n(k-1)-[\los+\rcs],\\
{\mak}'&=& \lob +\rcb,\\
{\lmak}'&=& n(k-1)-[\lcb+\rob].
\end{eqnarray*}
The following result was first noticed by Ksavrelof
and Zeng in \cite{KsZe}.
For completeness, we include a more straightforward proof.
\begin{prop} For any $\pi\in \OP_n^k$ we have
$$
\mak\,=\,\lmak{}'\quad\textrm{and}\quad \mak{}'\,=\,\lmak.
$$
\end{prop}
 \proof For $\pi=B_1-B_2-\cdots -B_k\in \OP_n^k$ and $i\in [n]$ we have
 $$
(\los_i+\lob_i+\ros_i+\rob_i)\pi=(\lcs_i+\lcb_i+\rcs_i+\rcb_i)\pi=k-1.
 $$
Il follows that
$$
\los+\lob+\ros+\rob=\lcs+\lcb+\rcs+\rcb=n(k-1).
$$
The Proposition is then equivalent  to
$$
\lob+\los=\lcb+\lcs,
$$
which is obvious.
\qed
\medskip

In view of the above proposition the conjectures in \cite{Stein} are reduced to the following
\begin{conj}[Steingr\'{\i}msson] The following
statistics are Euler-Mahonian on $\OP$ :
\begin{align*}
&\mak+\bInv,\quad\lmak+\bInv,\quad\mak+\bMaj,\quad\lmak+\bMaj,\\
&\cinvLSB:=\lsb + \cbInv + {k\choose2}  \quad and\quad
\cmajLSB:=\lsb + \cbmaj+ {k\choose2},
\end{align*}
where $\cbInv = {k\choose2}-\bInv$ and $\cbmaj ={k\choose2}-\bMaj$.
In other words, the  generating functions of the above statistics
over  $\OP_n^k$
are equal to  $[k]_q!S_q(n,k)$.
\end{conj}

\medskip
Let $\pi$ be a partition of $[n]$.
A \emph{singleton} is the
element of a block which has only one element. Now, consider
a  block $B$ of a partition $\pi$ whose cardinal is $\geq2$. An
element of $B$ is a {\it strict  opener} (resp. {\it strict closer})
if it is the least (resp. greatest) element of $B$, and a {\it
transient} if it is neither the least nor greatest element of $B$.

  The sets of strict openers, strict closers, singletons and  transients of $\pi$
will be denoted  by $\O(\pi)$, $\C(\pi)$, $\S(\pi)$ and $\T(\pi)$,
respectively. The $4$-tuple
$\lambda(\pi)=(\O(\pi),\T(\pi),\S(\pi),\C(\pi))$ is called
the \emph{type} of $\pi$.
  For instance, for the partition $\pi=3\,5-2\,4\,6-1-7\,8$, we get
  $$
  \O(\pi)=\{2,3,7\},\quad
\C(\pi)=\{5,6,8\}, \quad \S(\pi)=\{1\}\quad\textrm{and} \quad
\T(\pi)=\{4\}.
$$
  Clearly we have  $\open = \O\cup\S$ and $\clos = \C\cup\S$
therefore we get
\begin{align*}
&\bInv=\rcs(\O \cup \S) \quad\textrm{and}\quad
\inv=\ros(\O
 \cup\S),\\
 &\bExc=\lcs(\O \cup \S) \quad\textrm{and}\quad \exc=\los(\O
\cup\S).
\end{align*}

\subsection{Main results}
Consider the following two generating functions of ordered
partitions with $k\geq 0$ blocks:
\begin{align}
\phi_k(a;x,y,t,u):&=\sum_{\pi\,\in\,\OP^k}x^{(\mak+\bInv)\pi}\,y^{
\cinvLSB\,\pi}
\,t^{\inv\,\pi}\,u^{\exc\,\pi}a^{|\pi|},\label{eq:val1}\\
\varphi_k(a;z,t,u):&=\sum_{\pi\,\in\,\OP^k}z^{(\lmak+\bInv)\pi}\,t^{\inv\,\pi}
\,u^{\exc\,\pi}\,a^{|\pi|}.\label{eq:val2}
\end{align}
The following is the main result of this paper.
\begin{thm}\label{thm:rec}
We have
\begin{align}
\phi_k(a;x,y,t,u)=& \frac{a^{k}\,(xy)^{{k\choose2}}[k]_{tx,uy}!}
{\prod_{i=1}^k(1-a[i]_{x,y})},\label{eq:rec1}\\
\varphi_k(a;z,t,u)=&
\frac{a^kz^{k\choose 2}\,[k]_{tz,u}!}{\prod_{i=1}^k(1-a[i]_z)}\,.
\label{eq:rec2}
\end{align}
\end{thm}
The proof of this theorem will occupy the whole Section~3. We first
derive some results on Euler-Mahonian
statistics on ordered partitions.

By definition, the combinatorial interpretations of the
following specializations of $\phi_k$ and $\varphi_k$
 are obvious:
\begin{align*}
&\sum_{\pi\in\OP^k}q^{(\mak+\bInv)\pi}\,a^{|\pi|}=\phi_k(a;q,1,1,1),\\
&\sum_{\pi\in\OP^k}q^{\cinvLSB\,\pi}\,a^{|\pi|}=\phi_k(a;1,q,1,1),\\
&\sum_{\pi\in\OP^k}
q^{(\mak+\bInv-\inv+\exc)\pi}\,a^{|\pi|}=\phi_k(a;q,1,1/q,q),\\
&\sum_{\pi\in\OP^k}q^{(\cinvLSB+\inv-\exc)\,\pi}\,a^{|\pi|}=\phi_k(a;1,q,q,1/q),\\
&\sum_{\pi\in\OP^k}q^{(\lmak+\bInv)\pi}\,a^{|\pi|}=\varphi_k(a;q,1,1),\\
&\sum_{\pi\in\OP^k}q^{(\lmak+\bInv-\inv+\exc)\pi}\,a^{|\pi|}=\varphi_k(a;q,1/q,q).
\end{align*}
Applying Theorem~\ref{thm:rec} we see that
the right-hand sides of  the
above six identities are all equal to
\begin{equation}\label{eq:gfstirling}
\frac{a^{k}\,q^{{k\choose2}}[k]_q!}
{\prod_{i=1}^k(1-a[i]_q)}=\sum_{n\geq k}[k]_q!S_q(n,k)\,a^n,
\end{equation}
where the last equality follows directly
from \eqref{eq:stirling}. Thus we have proved

\begin{thm}\label{thm:cor1}
The following six inversion-like statistics are Euler-Mahonian on
$\OP$:
 \begin{align*}
  \mak +\bInv, \quad &\quad \mak + \bInv -(\inv-\exc)\;,\\
\lmak+\bInv,\quad &\quad \lmak +\bInv+(\inv-\exc)\;,\\
\cinvLSB,\quad &\quad \cinvLSB+(\inv-\exc)\;.
\end{align*}
\end{thm}
\subsection{Consequence on partitions}
Since a partition is an ordered partition without inversion, so we
can derive the following "hard" combinatorial interpretations for
$q$-Stirling numbers by putting $t= 0$ and extracting the coefficient of
$a^n$ in Theorem~\ref{thm:rec}:
$$
S_q(n,k)=\sum_{\pi\in\P_n^k}q^{\mak\,\pi}=\sum_{\pi\in\P_n^k}q^{\lmak\,\pi}=\sum_{\pi\in\P_n^k}q^{\lsb\,\pi+{k\choose
2}}.
$$
The first two interpretations were proved by Ksavrelof and Zeng
\cite{KsZe}.  The third
interpretation was first proved by Stanton (see \cite{WW}).

Note that by definition $k(1-k)+\cinvLSB=\lsb -\bInv$, then
by noticing that the two statistics $\inv$ and $\bInv$ vanish on
(unordered) partitions, we get that
\begin{align*}
&\sum_{\pi\in\P^k}q^{\mak\,\pi}a^{|\pi|}=\phi_k(a;q,1,0,1),\\
&\sum_{\pi\in\P^k}q^{{k\choose 2}+\lsb\,\pi}a^{|\pi|}=q^{k\choose2}\phi_k(a;1,q,0,1),\\
&\sum_{\pi\in\P^k}q^{\lmak\,\pi}a^{|\pi|}=\varphi_k(a;q,0,1).
\end{align*}
Now, applying Theorem~\ref{thm:rec} we see that
the right-hand-sides of  the
above three identities are equal to
$\sum_{n\geq k}S_q(n,k)\,a^n$ in view of
\eqref{eq:gfstirling}.

\section{Proof of Theorem 2.3}
\subsection{Ordered partitions and  walks in digraphs}

   Let $\pi=B_1-B_2-\cdots -B_k$ be a partition of $[n]$ and $i$ an integer in $[n]$.
The restriction $B_j(\leq i):=B_j \cap [i]$ of the block $B_j$ is
said to be \emph{opened} if $B\not\subseteq [i]$ and $B\cap
[i]\neq\emptyset$, \emph{closed} if $B\subseteq [i]$, and
\emph{empty} if $B\cap [i]=\emptyset$.
  The $i$-th \emph{trace} of $\pi$, $T_i(\pi)$, is defined by
$$
 T_i(\pi)=B_1(\leq i)-B_2(\leq i)-\cdots -B_k(\leq i),
$$
where the empty restrictions are not written. The sequence
$(T_i(\pi))_{1\leq i\leq n}$ is called the \emph{trace} of the
partition $\pi$.
 We denote by $\op_i\;\pi$ and $\cl_i\;\pi$ the numbers of opened blocks and closed
blocks, respectively, in $T_i(\pi)$ and set
$\form_i(\pi)=(\cl_i\,\pi
,\op_i \,\pi)$ for $1\leq i\leq n$ with $\form_0 (\pi) =(0,0)$. The
sequence $(\form_i(\pi))_{0\leq i\leq n}$ is called the \emph{form}
of the partition $\pi$.

 For instance, if {\small $\pi=\{6\}-\{3,5,7\}-\{1,4,10\}-\{9\}-\{2,8\}$},
then {\small
$T_6(\pi)=\{6\}-\{3,5,\cdots\}-\{1,4,\cdots\}-\{2,\cdots\}$}, where
each opened block has an ellipsis, and we get {\small
$\form_6(\pi)=(1,3)$}.\\

\begin{rmk}\label{rmk:form-type}
 Given the form of a partition, it is easy to deduce its type, and reversely.
Indeed, let $i\in[n]$ and suppose that $\form_{i-1}(\pi)=(k,l)$,
then
$$
  \form_i(\pi)=\left\{
\begin{array}{ll}
    (k,l+1), & \hbox{if $i\in \O(\pi)$;} \\
    (k+1,l), & \hbox{if $i\in \S(\pi)$;} \\
   (k,l), & \hbox{if $i\in \T(\pi)$;} \\
    (k+1,l-1), & \hbox{ if $i\in \C(\pi)$.} \\
\end{array}%
\right.
 $$
\end{rmk}
Note that if $l=0$, then $i$ can be neither a strict closer nor a
transient ($l=0$ means that all the blocks in $T_{i-1}\,\pi$ are
closed).

For any integer $k\geq0$,
let $D_k$  be the digraph with
 vertex set
  $V_k=\{(i,j)\in \mathbb{N}^2|\,i+j\leq k\}$,
  and there is an edge in
$D_k$ from $(i,j)$ to $(i',j')$ if and only if $(i',j')=(i,j)$ with
$j>0$ or $(i',j')\in\{(i,j+1),\;(i+1,j),\;(i+1,j-1)\}$.
 It is obvious that the number of vertices of $D_k$ is equal to
$$
\N{k}:=1+2+\cdots+(k+1)=\frac{(k+1)(k+2)}{2}.
$$
Let $v_1,\cdots,v_{\N{k}}$ be
the vertices of $D_k$
arranged according to the following order:
$(i,j)\leq(i',j')$ if and only if
$i+j<i'+j'$ or ($i+j=i'+j'$ and $j\geq j'$).
For
instance, we get
$v_1=(0,0),\,v_2=(0,1),\,v_3=(1,0),\,v_4=(0,2),\,v_5=(1,1),
\,v_6=(2,0),\cdots,v_{\N{k}}=(k,0)$.
An illustration of $D_k$ is given in Figure~2.

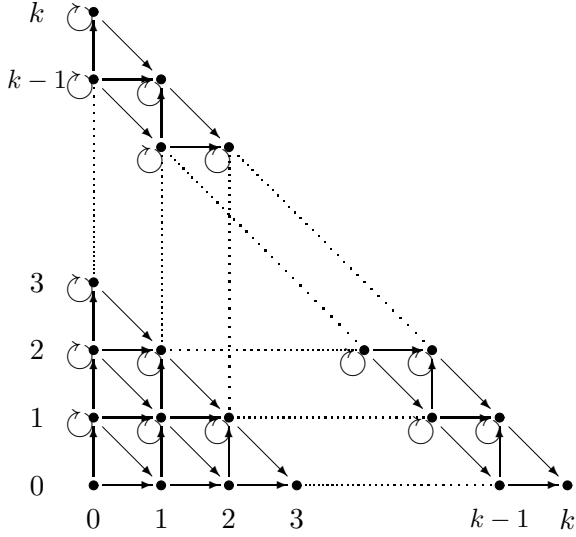
\begin{figure}[h]
\begin{center}
{\setlength{\unitlength}{0.60mm}
\begin{picture}(130,120)(0,-10)
\put(0,0){\circle*{2}}\put(15,0){\circle*{2}}\put(30,0){\circle*{2}}\put(45,0){\circle*{2}}
\put(0,15){\circle*{2}}\put(15,15){\circle*{2}}\put(30,15){\circle*{2}}
\put(0,30){\circle*{2}}\put(15,30){\circle*{2}}\put(0,45){\circle*{2}}
\put(0,105){\circle*{2}}\put(0,90){\circle*{2}}\put(15,90){\circle*{2}}\put(15,75){\circle*{2}}\put(30,75){\circle*{2}}
\put(105,0){\circle*{2}}\put(90,0){\circle*{2}}\put(90,15){\circle*{2}}\put(75,15){\circle*{2}}\put(75,30){\circle*{2}}
\put(60,30){\circle*{2}}

\put(2,0){\vector(1,0){11}}\put(17,0){\vector(1,0){11}}\put(32,0){\vector(1,0){11}}\put(62,30){\vector(1,0){11}}
\put(0,2){\vector(0,1){11}}\put(0,17){\vector(0,1){11}}\put(0,32){\vector(0,1){11}}
\put(2,13){\vector(1,-1){11}}\put(2,28){\vector(1,-1){11}}\put(2,43){\vector(1,-1){11}}
\put(17,13){\vector(1,-1){11}}\put(17,28){\vector(1,-1){11}}\put(32,13){\vector(1,-1){11}}
\put(2,15){\vector(1,0){11}}\put(17,15){\vector(1,0){11}}\put(2,30){\vector(1,0){11}}
\put(92,0){\vector(1,0){11}}\put(77,15){\vector(1,0){11}}\put(2,90){\vector(1,0){11}}\put(17,75){\vector(1,0){11}}
\put(0,92){\vector(0,1){11}}\put(15,77){\vector(0,1){11}}\put(15,2){\vector(0,1){11}}\put(15,17){\vector(0,1){11}}
\put(30,2){\vector(0,1){11}}\put(90,2){\vector(0,1){11}}\put(75,17){\vector(0,1){11}}
\put(2,103){\vector(1,-1){11}}\put(2,88){\vector(1,-1){11}}\put(17,88){\vector(1,-1){11}}
\put(62,28){\vector(1,-1){11}}\put(77,28){\vector(1,-1){11}}\put(77,13){\vector(1,-1){11}}\put(92,13){\vector(1,-1){11}}
\bezier{30}(17,73)(30,60)(58,32)\bezier{30}(0,47)(0,60)(0,88)\bezier{30}(47,0)(60,0)(88,0)\bezier{30}(32,73)(47,58)(73,32)
\bezier{30}(32,15)(50,15)(73,15)\bezier{30}(17,30)(47,30)(58,30)\bezier{30}(15,32)(15,50)(15,73)\bezier{30}(30,17)(30,53)(30,73)
\put(0,15){\makebox(-6,-3){\large$\circlearrowright$}}\put(0,30){\makebox(-6,-3){\large$\circlearrowright$}}
\put(0,45){\makebox(-6,-3){\large$\circlearrowright$}}\put(0,90){\makebox(-6,-3){\large$\circlearrowright$}}
\put(0,105){\makebox(-6,-3){\large$\circlearrowright$}}\put(15,15){\makebox(-5,-6){\large$\circlearrowright$}}
\put(15,30){\makebox(-5,-6){\large$\circlearrowright$}}\put(15,75){\makebox(-5,-6){\large$\circlearrowright$}}
\put(15,90){\makebox(-5,-6){\large$\circlearrowright$}}\put(30,15){\makebox(-5,-6){\large$\circlearrowright$}}
\put(30,75){\makebox(-5,-6){\large$\circlearrowright$}}\put(60,30){\makebox(-5,-6){\large$\circlearrowright$}}
\put(75,30){\makebox(-5,-6){\large$\circlearrowright$}}\put(75,15){\makebox(-5,-6){\large$\circlearrowright$}}
\put(90,15){\makebox(-5,-6){\large$\circlearrowright$}}
\put(0,0){\makebox(0,-15){\small$0$}}\put(15,0){\makebox(0,-15){\small$1$}}\put(30,0){\makebox(0,-15){\small$2$}}
\put(45,0){\makebox(0,-15){\small$3$}}\put(90,0){\makebox(0,-15){\footnotesize$k-1$}}\put(105,0){\makebox(0,-15){\small$k$}}
\put(0,0){\makebox(-25,0){\small$0$}}\put(0,15){\makebox(-25,0){\small$1$}}\put(0,30){\makebox(-25,0){\small$2$}}
\put(0,45){\makebox(-25,0){\small$3$}}\put(0,90){\makebox(-25,0){\footnotesize$k-1$}}\put(0,105){\makebox(-25,0){\small$k$}}
\end{picture}}
\end{center}
\caption{The digraph $D_k$}
\end{figure}

 A path of length $n$ in $D_k$ is a finite
sequence $(s_0, s_1, \ldots, s_n)$ of points in $V_k$ such that
$(s_i, s_{i+1})$ is an edge for $i=0,\ldots,n-1$.
The step $(s_i,s_{i+1})$ is called \emph{ North} (resp. \emph{East},
 \emph{South-East} and \emph{Null})
  if $s_{i+1}=(x_i,y_i+1)$
 (resp. $s_{i+1}=(x_i+1,y_i)$, $s_{i+1}=(x_i+1,y_i-1)$ and $s_{i+1}=s_i$).
 The number $x_i$ and $y_i$ are respectively the \emph{abscissa} and the \emph{height} of the
 step $(s_i,s_{i+1})$.

\begin{defn} A path in $D_k$ from
$v_1=(0,0)$ to $v_{\N{k}}=(k,0)$ is called a path of depth $k$.
 Let $\Omega_n^k$ be the set of paths of depth $k$ and length $n$ and
$\Omega^k=\bigcup_{n\geq 0}\Omega_n^k$ be the set of all paths of
depth $k$.
\end{defn}
 The following result just follows from the definition of paths.
\begin{prop}\label{prop:form-path}
For $n\geq k\geq 0$, the forms of partitions in $\OP_n^k$ are
exactly the paths in $\Omega_n^k$.
\end{prop}

We can visualize a path by drawing a segment from
$s_{i-1}$ to $s_i$ in the x-y plan. For
instance,
 the path
$$
\omega=((0,0),\,(0,1),\,(0,2),\,(0,3),\,(0,3),\,
(0,3),\,(1,3),\,(2,2),\,(3,1),\,(4,1),\,(5,0))
$$
is illustrated in Figure~2.
\begin{figure}[h]\
\begin{center}
{\setlength{\unitlength}{0.4mm}
\begin{picture}(90,60)(0,0)
\put(0,0){\circle*{3}}\put(0,15){\circle*{3}}\put(0,30){\circle*{3}}\put(0,45){\circle*{3}}\put(15,45){\circle*{3}}
\put(30,30){\circle*{3}}\put(45,15){\circle*{3}}\put(60,15){\circle*{3}}\put(75,0){\circle*{3}}
\put(-10,0){\vector(1,0){100}}\put(0,-10){\vector(0,1){75}}\put(0,0){\linethickness{0.2mm}\line(0,1){45}}\put(0,45){\linethickness{0.2mm}\line(1,0){15}}
\put(15,45){\linethickness{0.2mm}\line(1,-1){30}}\put(45,15){\linethickness{0.2mm}\line(1,0){15}}\put(60,15){\linethickness{0.2mm}\line(1,-1){15}}
\put(-5,45){\linethickness{0.2mm}\circle{10}}\put(-15,46){\scriptsize2}
\put(-4,-7){\scriptsize 0}\put(74,-7){\scriptsize5}
\end{picture}}
\end{center}
\caption{A path in $\Omega_{10}^5$, where $2$ means two successive Null
steps at $(0,3)$. }
\end{figure}
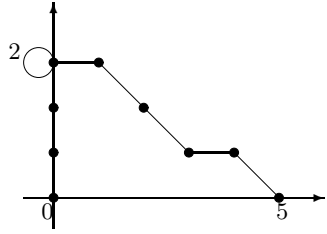

   By definition, an ordered partition on $[i]$ with two kinds of blocks, opened or closed,
is called a \emph{trace} on  $[i]$.
  Let $T_{i-1}$ be a trace on $[i-1]$ and suppose that $\form\,(T_{i-1})= (k,l)$,
with $k,l\geq 0$. There are several possibilities to insert the
element $i$ in $T_{i-1}$ according to the nature of $i$. The element
$i$ could be~:
\begin{itemize}
\item[(a)]a strict closer (resp. a transient): there are $p_i=l$ possibilities
 to close with $i$ (resp. insert $i$ in) one of the $l$ opened blocks of $T_{i-1}$.
\item[(b)]a singleton (resp. a strict opener): there are  $p_i=k+l+1$ possibilities to insert the
singleton $\{i\}$ (resp. open a block with $i$) in $T_{i-1}$ (before
all the blocks, between two blocks, or after all the blocks).
\end{itemize}
We observe that if $l=0$, then in case $(a)$, there are none
possibility to add $i$ ($p_i=0$), which is natural because all the
blocks in $T_{i-1}$ are closed.

\begin{rmk}\label{rmk:choice}
  The number of ways $p_i$ to add the element $i$, according to its "nature",
 in a trace $T_{i-1}$ on $[i-1]$,
depends only on $\form\,(T_{i-1})$. If $\form\,(T_{i-1})=(k,l)$,
then $p_i=l$ (resp. $p_i= k+l+1$) if we want insert $i$ as a
transient or a strict closer (resp. a strict opener or a singleton).
\end{rmk}

 We assume the possibilities to add an element (according its
nature) in a trace are arranged from left to right. Namely, if we
insert a singleton or open a block (resp. insert a transient or
close a block), the spaces (resp. opened blocks) which correspond to
the
possibilities are arranged from left to right.\\

 For instance, $T_6=\{6\}-\{3,5,\cdots\}-\{1,4,\cdots\}-\{2,\cdots\}$
is a trace on $[6]$, then $\form\,(T_6)$=(1,3) and there are:
\begin{itemize}
\item[(i)] $5$ possibilities to open a block or insert a singleton in $T_6$.
$$
\begin{array}{cccccccccc}
T_5=&        &\{6\}&     &\{3,5,\cdots\}&  &\{1,4,\cdots\}&   &\{2,\cdots\}&\\
    &\uparrow&     &\uparrow&       &\uparrow&    &\uparrow& &\uparrow\\
$\footnotesize{ choice}$&1     &     &  2    &          &3   &      &4      &  &5\\
\end{array}
$$

\item[(ii)] $3$ possibilities to close a block or add a transient
in $T_6$. Namely,
$$
\begin{array}{ccccc}
T_6=&\{6\}&\{3,5\cdots\}&\{1,4,\cdots\}&\{2,\cdots\}\\
&&\uparrow&\uparrow&\uparrow\\
$\footnotesize{ choice}$&&1&2&3
\end{array}
$$
\end{itemize}

\begin{defn}
 A path diagram of depth $k$ and length $n$ is a pair
$(\omega,\xi)$, where $\omega$ is a path of
depth $k$ and length $n$ and $\xi=(\xi_i)_{1\leq i\leq n}$ is a
sequence of integers such that $1\leq \xi_i\leq q$ if the $i$-th
step of $\omega$ is Null or South-East of height $q$, and $1\leq
\xi_i\leq p+q+1$ if the $i$-th step of $\omega$ is North or East
of abscissa $p$ and height $q$.
\end{defn}
 Denote by $\Delta_n^k$ the set of path diagrams of depth $k$ and length $n$
and by $\Delta^k=\bigcup_{n\geq 0}\Delta_n^k$ the set of path
diagrams of depth $k$.\\

{\bf The mapping $\psi$:}
 Let $n\geq k\geq1$. Given a path
 diagram $h=(\omega,\xi)\in\Delta_n^k$, we associate a partition $\psi(h)\in\OP_n^k$ by
 constructing successively its traces $T_i$ for $1\leq i\leq n$ as follows:
\begin{enumerate}
\item Set $T_0=\emptyset$.
\item For $1\leq i\leq n$, we construct $T_i$ from $T_{i-1}$ by the following process.
 Suppose that $s_{i-1}=(p,q)$ and the $i$-th step of $\omega$ is :
 \begin{itemize}
\item [(i)] North (resp. East), then we open a block with $i$
(resp. insert the singleton $\{i\}$) in $T_{i-1}$ according to the
choice $\xi_i$.
\item[(ii)] South-East (resp. Null), then we close with $i$
(resp. insert $i$ as a transient in ) an opened block of $T_{i-1}$
according to the choice $\xi_i$.
  \end{itemize}
\item Set $\psi(h)=T_n$
\end{enumerate}
For instance, if $h=(\omega,\xi)$ where $\omega$ is the path of
Figure~1 and {\small $\xi=(1,2,1,2,1,1,1,2,4,1)$}, then
the step by step construction of $\psi(h)$ goes as follows:
$$
\scriptsize{
\begin{array}{ccccccl}
i && \quad\textrm{step}_i\quad& \xi_i &&& \quad T_i\\
&&&&&&\\
1 && North& 1 &&& \{1,\cdots\} \\
&&&&&&\\
2 && North& 2 &&& \{1,\cdots\}-\{2,\cdots\} \\
&&&&&&\\
3 && North& 1 &&& \{3,\cdots\}-\{1,\cdots\} -\{2,\cdots\} \\
&&&&&&\\
4 && Null & 2 &&& \{3,\cdots\}-\{1,4,\cdots\}-\{2,\cdots\} \\
&&&&&&\\
5 && Null & 1 &&& \{3,5,\cdots\}-\{1,4,\cdots\}-\{2,\cdots\} \\
&&&&&&\\
6 && East & 1 &&& \{6\}-\{3,5,\cdots\}-\{1,4,\cdots\}-\{2,\cdots\} \\
&&&&&&\\
7 && South\text{-}East& 1 &&& \{6\}-\{3,5,7\}-\{1,4,\cdots\}-\{2,\cdots\} \\
&&&&&&\\
8 && South\text{-}East& 2 &&& \{6\}-\{3,5,7\}-\{1,4,\cdots\}-\{2,8\} \\
&&&&&&\\
9 && East & 4 &&& \{6\}-\{3,5,7\}-\{1,4,\cdots\}-\{9\}-\{2,8\}\\
&&&&&&\\
10 && South\text{-}East & 1
&&&\{6\}-\{3,5,7\}-\{1,4,10\}-\{9\}-\{2,8\}.
\end{array}}
$$
Thus
$\psi(h)= \{6\}-\{3,5,7\}-\{1,4,10\}-\{9\}-\{2,8\}$.
\begin{thm}\label{thm:bij}
 For each $n\geq k\geq1$, the mapping
 $\psi: h=(\omega,\xi)\mapsto \pi$  is a bijection from
 $\Delta_n^k$ to $\OP_n^k$ such that:\\
  if the $i$-th
 step of $\omega$
is North or East of abscissa $p$ and height $q$, then $i\in (\O\cup\S)(\pi)$ and
$$(\lcs+\rcs)_i(\pi)=p,\quad  (\lsb+\rsb)_i(\pi)=q,\quad
    \los_i(\pi)=\xi_i-1 \quad \text{and}\quad
    \ros_i(\pi)=p+q+1-\xi_i;
    $$
if the $i$-th step of $\omega$ is  South-East or Null  of abscissa $p$ and height $q$, then
$i\in (\T\cup\C)(\pi)$ and
$$
(\lcs+\rcs)_i(\pi)=p,\quad (\lsb+\rsb)_i(\pi)=q-1,\quad
\lsb_i(\pi)=\xi_i-1 \quad \text{and}\quad \rsb_i(\pi)=q-\xi_i.
$$
\end{thm}
\pf
By Remark \ref{rmk:choice} and Proposition~\ref{prop:form-path}, it is easy
to see that the above algorithm is well defined.
Suppose we are constructing the ordered partition $\pi$ and we
arrive at the $i$-th step of the construction. The $i$-th step of
$w$ is a step with initial vertex $(p,q)$. At this step of the
construction, there are exactly $p+q$ blocks in $T_{i-1}$, whose $p$
(resp. $q$) are closed (resp. opened), and all the elements in
$T_{i-1}$ are strictly inferior than $i$. Suppose that the $i$-th
step of $w$ is a step of type:
\begin{itemize}
\item[(i)] North or East : then $1\leq a_i\leq p+q+1$ and
 by remark \ref{rmk:form-type}, the element
$i$ will be a strict opener or a singleton of the partition $\pi$.
It's clear that $(\lcs+\rcs)_i\,(\pi)$ (resp. $(\lsb+\rsb)_i\,(\pi$)) is
equal to the number of closed (resp. opened) blocks in
$T_{i-1}\,(\pi)$. Thus, $(\lcs+\rcs)_i\,(\pi)=p$ and
$(\lsb+\rsb)_i\,(\pi)=q$.
\\ Because all the blocks in $T_{i-1}$ have opener smaller than i,
then $\los_i\,\pi$ (resp. $\ros_i\,\pi$) is just the number of blocks
in $T_i(\pi)$ on the left (resp. right) of the block which contains
$i$. Thus, because we open a block with the element $i$ or add the
singleton $\{i\}$ between the $(\xi_i-1)$-th and $(\xi_i)$-th blocks
of $T_{i-1}$, we get $(\los_i,\ros_i)(\pi)=(\xi_i-1,p+q+1-\xi_i)$ .

\item[(ii)] South-East or Null : then $1\leq a_i\leq q$,
 and by remark \ref{rmk:form-type}, the element $i$ will be a transient or a
strict closer of the partition $\pi$. By the same arguments as case
$(i)$, we get $(\lcs+\rcs)_i\,(\pi)=p$. Remark $(\lsb+\rsb)_i\,(\pi)$ is
equal to the number of opened blocks in $T_{i}\,(\pi)$ which don't
contain the element $i$, thus because we insert $i$ in one of the
opened block in $T_{i-1}$, we get $(\lsb+\rsb)_i\,(\pi)=q-1$.
Moreover, $\lsb_i\,(\pi)$ (resp. $\rsb_i\,(\pi)$) is equal to the
number of opened blocks in $T_i(\pi)$ on the left (resp. on the
right) of the block which contains $i$. Because we insert $i$ in the
$\xi_i$-th opened block in $T_{i-1}(\pi)$, we get
$(\lsb_i,\rsb_i)(\pi)=(\xi_i-1,q-\xi_i)$.
\end{itemize}
\qed

\subsection{Generating functions of walks}
 For $0\leq k\leq n$, let $\textbf{t}=(t_1,t_2,t_3,t_4,t_5,t_6,t_7)$
 and
\begin{align}
Q_{n,k}(\textbf{t}):&=
\sum_{\pi\in\OP_n^k}t_1^{(\lcs+\rcs)(\O\cup\S)\pi}\,
t_2^{(\lcs+\rcs)(\T\cup\C)\pi}\,t_3^{\rsb(\T \cup
\C)\pi}\label{eq:defQ}\\
&\hspace{2cm}\times \,t_4^{\lsb(\T \cup \C)\pi}\,t_5^{\ros(\O \cup
\S)\pi}\,t_6^{\los(\O \cup
\S)\pi}\,t_7^{(\lsb+\rsb)(\O\cup\S)\pi}.\nonumber
\end{align}

 Given a path $\omega$, define the weight $\text{v}(\omega)$ of $\omega$
to be the product of the weights of all its steps, where the weight
of a step of abscissa $i$ and height $j$ is:
\begin{equation}\label{eq:val}
v(\omega)=\left\{
  \begin{array}{ll}
    t_1^i\,t_7^j\,[i+j+1]_{t_5,t_6} & \hbox{if the step is North or East;} \\
    t_2^i\,[j]_{t_3,t_4} & \hbox{if the step is Null or South-East.}
  \end{array}
\right.
\end{equation}

It follows easily from Theorem~\ref{thm:bij} that
\[
\sum_{\omega\in\Omega_{n}^k}\text{v}(\omega)=Q_{n,k}(\textbf{t}) \;.
\]

 Denote by $|\omega|$ the length of the path $\omega$.
Then, using the above identity, we get
  \begin{equation}\label{eq:codage}
  Q_k(a;\textbf{t}):=\sum_{n\geq
0}Q_{n,k}(\textbf{t})\,a^n=\sum_{w\in\Omega^k}\text{v}(w)\,a^{|\omega|}\;.
  \end{equation}

 The \emph{adjacency matrix} $A_k$ of $D_k$ relative
to the valuation $\text{v}$ is the $\N{k}\times \N{k}$ matrix
defined by
$$A_k(i,j)=\left\{
            \begin{array}{ll}
              \textrm{v}(v_i,v_j) & \hbox{if $(v_i,v_j)$ is an edge of $D_k$;} \\
              0 & \hbox{otherwise.}
            \end{array}
          \right.
$$
Applying transfer-matrix method (see e.g. \cite[Theorem 4.7.2]{Stan}), we derive
\begin{equation}\label{eq:transfert}
Q_k(a;\textbf{t})=\frac{(-1)^{1+\N{k}}\det (I-aA_k;
\N{k},1)}{\det(I-aA_k)},
\end{equation}
where $(B;i,j)$ denotes the matrix obtained by removing the $i$-th
row and $j$-th column of $B$ and $I$ is the $\N{k}\times \N{k}$
identity matrix.

In order to
prove Theorem \ref{thm:rec} we need to evaluate \eqref{eq:transfert}
in the following special cases:
\begin{align}
f_k(a;x,y,t,u)=Q_k(a;x,x,x,y,t,u,y),\label{eq:spe1}\\
g_k(a;z,t,u)=Q_k(a;1,z,1,z,t,u,1).\label{eq:spe2}
\end{align}

  Let $A_k'$ and  $A_k''$ be the adjacency matrix of $D_k$ relative to the
weight function $\textrm{v}'$ and   $\textrm{v}''$
obtained from the weight function $\textrm{v}$ by making the
substitution \eqref{eq:spe1} and \eqref{eq:spe2}, respectively.
 Namely, the weights $\textrm{v}'(e)$ and $\textrm{v}''(e)$
 of an edge $e=((i,j), (i',j'))$ of $D_k$
with initial vertex $(i,j)$ are~:
$$
\textrm{v}'(e)=\left\{
  \begin{array}{ll}
    x^i\,y^j\,[i+j+1]_{t,u} & \hbox{if $(i',j')=(i,j+1)$ or
$(i+1,j)$ ;} \\
    x^i\,[j]_{x,y} & \hbox{if $(i',j')=(i,j)$ or $(i+1,j-1)$,}
  \end{array}
\right.
$$
and
$$\textrm{v}''(e)=\left\{
  \begin{array}{ll}
    [i+j+1]_{t,u} & \hbox{if $(i',j')=(i,j+1)$ or
$(i+1,j)$ ;} \\
    z^i[j]_z & \hbox{if $(i',j')=(i,j)$ or $(i+1,j-1)$.}
  \end{array}
\right.
$$

Now, for each $k\geq 0$ let
$$
M_k=I-aA_k' \quad \textrm{and}\quad N_k=I-aA_k''.
$$
 Then by \eqref{eq:transfert}, \eqref{eq:spe1} and \eqref{eq:spe2}
  we have
\begin{align}
f_{k}(a;x,y,t,u) & =\frac{(-1)^{1+\N{k}}\det(M_k;
\N{k},1)}{\det M_k}, \label{eq:transfertspe1}\\
g_k(a;z,t,u)
&=\frac{(-1)^{1+\N{k}}\det (N_k; \N{k},1)}{\det N_k}
\label{eq:transfertspe2}
\end{align}
for each $k\geq 1$.

Since $M_n$ and $N_n$ are upper triangular matrices it is easy to
see that for each $n\geq1$
\begin{align}
&\det M_n =\prod_{m=1}^{n}\prod_{i=0}^{m}(1-ax^{i}[m-i]_{x,y}),\\
&\det N_n=\prod_{m=1}^{n}\prod_{k=0}^{n-m}(1-az^{k}[m]_{q}).
\end{align}

The evaluation of $\det(M_n;\N{n},1)$ and $\det(N_n;\N{n},1)$
is not simple(see last section).
\begin{thm}\label{thm:minor1}
Let $n\geq1$ be a positive integer. Then
\begin{align}
\det(M_n;\N{n},1)
&=(-1)^{n\choose2}a^nx^{n\choose2}[{n}]_{t,u}!\prod_{m=1}^{n-1}
\prod_{i=1}^{m}(1-ax^{i}[m-i+1]_{x,y}), \label{eq:minor1}
\end{align}
\end{thm}
\begin{thm}\label{thm:minor2}
Let $n\geq1$ be a positive integer. Then
\begin{align}
\det (N_n; \N{n},1) &=(-1)^{n\choose2}a^{n}\,[n]_{t,u}!
\prod_{m=1}^{n-1}\prod_{k=1}^{n-m}(1-az^{k-1}[m]_{z}).\label{eq:minor2}
\end{align}
\end{thm}

It is now trivial to obtain the following result.
\begin{cor}\label{thm:pre}
  For $k\geq 0$, we have
\begin{align}
f_{k}(a;x,y,t,u) &= \frac{a^{k}x^{{k\choose2}}[k]_{t,u}!}
{\prod_{i=1}^k(1-a[i]_{x,y})},\label{eq:pre1}\\
g_{k}(a;z,t,u)&=\frac{a^k\,[k]_{t,u}!}{\prod_{i=1}^k(1-az^{k-i}[i]_z)}\,.
\label{eq:pre2}
\end{align}
\end{cor}

To see the relation between
 $f_k$ and $\phi_k$, $g_k$ and $\varphi_k$ we need to
 establish the following
\begin{lem}\label{lem:reecriture}
The following functional identities hold on $\OP_n^k$:
\begin{align*}
\mak+\bInv&=(\lcs +\rcs )+  \rsb (\T \cup \C )+\inv, \\
\lmak+\bInv&=n(k-1)-(\lcs +\rcs )(\T \cup \C)-\lsb(\T \cup \C)-\exc,\\
\cinvLSB&=(\lsb + \rsb)(\O \cup \S)+\lsb (\T\cup\C)+\inv+2\,\exc.
\end{align*}
\end{lem}
\begin{proof} By definition we have
   \begin{align*}
\mak+\bInv &=\lcs +\ros + \rcs(\O \cup \S)\\
&=(\lcs +\rcs )+\ros +(\rcs (\O \cup \S)-\rcs)\\
           &= (\lcs +\rcs )+ \ros -\rcs (\T \cup \C)\\
           &= (\lcs +\rcs )+  (\ros-\rcs ) (\T \cup \C )+\ros (\O \cup \S)\\
           &= (\lcs +\rcs )+  \rsb (\T \cup \C )+\ros (\O \cup\S).
           \end{align*}
Also
           \begin{align*}
n(k-1)-(\lmak+\bInv) &=(\los +\rcs ) - \rcs(\O \cup \S) = \los  +\rcs(\T \cup \C )\\
                   &= \los (\O \cup \S) +\rcs (\T \cup \C)+\los(\T \cup\C)\\
                   &= \exc+ (\lcs +\rcs )(\T \cup \C )+ (\los -\lcs)(\T \cup \C)\\
                   &= \exc+(\lcs +\rcs )(\T \cup \C )+ \lsb(\T \cup\C),
\end{align*}
and
\begin{align*}
\cinvLSB&=k(k-1)+\lsb-\bInv\\
&=2(\inv+\exc)+\lsb-\rcs (\O \cup \S)\\
         & =2(\inv+\exc)+ (\lsb + \rsb)(\O \cup \S)\\
         &\qquad \qquad-(\rcs+\rsb)(\O \cup \S)+\lsb (\T \cup \C)\\
         &= 2(\inv+\exc)+(\lsb + \rsb)(\O \cup \S)-\ros (\O \cup\S)+\lsb(\T \cup \C)\\
         &= (\lsb + \rsb)(\O \cup \S)+\lsb(\T \cup \C)+\inv+2\,\exc.
\end{align*}
The proof is thus completed.
\end{proof}

Now, we derive from \eqref{eq:spe1} and \eqref{eq:spe2} that
\begin{align*}
f_k(a;x,y,t,u)&=\sum_{\pi\in\OP^k}
x^{(\lcs+\rcs+\rsb(\T\cup\C))\,\pi}y^{((\lsb+\rsb)(\O\cup\S)+\lsb(\T\cup\C))\,\pi}
t^{\inv\,\pi}u^{\exc\,\pi}a^{|\pi|},\\
g_k(a;z,t,u)&=\sum_{\pi\in\OP^k}z^{(\lcs+\rcs+\lsb)(\T\cup\C)\pi}
t^{\inv\,\pi}u^{\exc\,\pi}a^{|\pi|}.
\end{align*}
It follows from the above lemma the following
\begin{lemma}\label{lem:substitution}
 The following identities hold:
\begin{align}
\phi_k(a;x,y,t,u)&=f_k(a;x,y,xyt,uy^2),\\
\varphi_k(a;z,t,u)&=g_k(az^{k-1};1/z,t,u/z).
\end{align}
\end{lemma}
Finally  Theorem \ref{thm:rec} follows immediately from Corollary~\ref{thm:pre}
and Lemma~\ref{lem:substitution}.
Therefore in order to prove Theorem~\ref{thm:rec} it remains to
prove Theorem~\ref{thm:minor1} and  Theorem~\ref{thm:minor2}.

\subsection{Proof of  Theorem \ref{thm:minor1}}
The matrix $M_n$ can be defined recursively by
\begin{align}
M_{0}=\left( 1\right),\quad
M_{n}=\left(
\begin{array}{c|c}
M_{{n}-1} &  \UM_{{n}-1} \\
\noalign{\medskip} \hline \noalign{\medskip} O_{{n}+1,\N{{n}-1}} &
\LM_{{n}-1}
\end{array}
\right), \label{eq:recM}
\end{align}
where $n\geq 1$,
\begin{equation}
\LM_{n-1}=\left(
\delta_{ij}-ax^{i-1}[n+1-i]_{x,y}(\delta_{ij}+\delta_{i+1,j})
\right)_{1\leq i,j\leq {n}+1}
\end{equation}
and $\UM_{{n}-1}$ is the $\N{{n}-1}\times ({n}+1)$ matrix
\begin{equation*}
\UM_{n-1}=\left(\begin{array}{c}
O_{\N{{n}-2},n+1}\\
\noalign{\medskip} \hline \noalign{\medskip} \UUM_{{n}-1}
\end{array}\right)
\end{equation*}
with the ${n}\times({n}+1)$ matrix
\begin{equation}
\UUM_{{n}-1}= \left(
-ax^{i-1}y^{{n}-i}[n]_{t,u}(\delta_{ij}+\delta_{i+1,j})
\right)_{1\leq i\leq {n},\,1\leq j\leq{n}+1}. \label{eq:UM}
\end{equation}
Here $\delta_{ij}$ stands for the Kronecker delta and $O_{m,n}$
denotes the $m\times n$ zero matrix.
%
%
For instance, we get
\[
M_{1}=\left( \scriptsize{
\begin{array}{c|cc}
1&-a&-a\\\hline
&&\\
0&1-a&-a\\
&&\\
0&0&1\\
\end{array}
} \right)
\]
and
\[
M_{2}=\left( \scriptsize{
\begin{array}{ccc|ccc}
1&-a&-a&0&0&0\\
&&&&&\\
0&1-a&-a&-ay(t+u)&-ay(t+u)&0\\
&&&&&\\
0&0&1&0&-ax(t+u)&-ax(t+u)\\
&&&&&\\
\hline
&&&&&\\
0&0&0&1-a(x+y)&-a(x+y)&0\\
&&&&&\\
0&0&0&0&1-ax&-ax\\
&&&&&\\
0&0&0&0&0&1
\end{array}
} \right).
\]

Let
\[
\K{n}=\N{n}-1=\frac{{n}({n}+3)}{2},
\]
and let $P_n=(M_n;\N{n},1)$, i.e the $\K{n}\times\K{n}$ matrix
obtained from $M_{n}$ by deleting the $\N{n}$th row and the first
column. $P_n$ can be defined as follows.
\[
P_{n}= \left(\begin{array}{c|c} P_{{n}-1}&\Pa_{{n}-1}\\
\noalign{\medskip} \hline \noalign{\medskip} X_{{n}-1}&\Pb_{{n}-1}
\end{array}
\right)
\]
 Here $\Pa_{{n}-1}$ is a $\K{{n}-1}\times({n}+1)$ matrix,
$X_{{n}-1}$ is a $({n}+1)\times \K{{n}-1}$ matrix, and $\Pb_{{n}-1}$
is a $({n}+1)\times({n}+1)$ matrix.  We shall compute $\det P_{n}$
by the following well-known formula for any block matrix with an
invertible square matrix $A$,
\[
\det\left(\begin{array}{c|c} A&B\\\noalign{\medskip}\hline
\noalign{\medskip} C&D
\end{array}\right)
= \det A\cdot \det\left(D-CA^{-1}B\right).
\]
Since the entries of $CA^{-1}B$ are also written by minors, we guess
these entries and prove it by induction (see Theorem~\ref{th:main1}).
Before we proceed to the proof, we give some examples.
\[
P_{1}=\left(
\begin{array}{cc}
-a&-a\\
&\\
1-a&-a\\
\end{array}
\right)
\]
and
\[
P_{2}=\left( \scriptsize{
\begin{array}{cc|ccc}
-a&-a&0&0&0\\
&&&&\\
1-a&-a&-ay(t+u)&-ay(t+u)&0\\
&&&&\\
\hline
&&&&\\
0&1&0&-ax(t+u)&-ax(t+u)\\
&&&&\\
0&0&1-a(x+y)&-a(x+y)&0\\
&&&&\\
0&0&0&1-ax&-ax\\
\end{array}
} \right).
\]
Thus, looking at $P_{2}$ as the block matrix composed of  $P_{1}$, $X_{1}$, $\Pa_{1}$
and $\Pb_{1}$, we have
\[
\scriptsize{ \Pa_{1} = \left( \begin {array}{ccc}
0&0&0\\
&&\\
-ay (t+u) &-ay(t+u) &0\\
\end {array} \right)
}
\]
and
\[
\scriptsize{\Pb_{1} = \left( \begin {array}{ccc}
0&-ax(t+u) &-ax(t+u) \\
&&\\
1-a ( x+y ) &-a ( x+y ) &0\\
&&\\
0&1-ax&-ax
\end {array} \right)
}.
\]
Since $\Pa_{n}$ is an $\K{n}\times ({n}+2)$ matrix, we can write
$$
\Pa_{n}= \left(\begin{array}{c} O_{\K{{n}-1},{n}+2}\\\hline U_{n}
\end{array}\right),
$$
where $U_{n}$ is the $({n}+1)\times({n}+2)$ matrix composed of the
last $({n}+1)$ rows of $\Pa_{n}$. For $1\leq {k}\leq {n}+2$, let
\[
P_{n}^{k} =\left(\begin{array}{c|c} P_{{n}-1}&\Pa_{{n}-1}\\
\noalign{\medskip} \hline \noalign{\medskip}
X_{{n}-1}&\Pb_{{n}-1}^{k}
\end{array}\right)
\]
denote the $\K{n}\times\K{n}$ matrix obtained from $P_{n}$ by
replacing the right-most column with the ${k}$th column of
$\Pa_{n}$. Here $\Pb_{{n}-1}^{k}$ is the $({n}+1)\times({n}+1)$
matrix obtained from $\Pb_{{n}-1}$ by replacing the right-most
column with the ${k}$th column of $U_{n}$. For example,
\[
 P_{2}^{2}=\left( \scriptsize{
\begin{array}{cc|ccc}
-a&-a&0&0&0\\
&&&&\\
1-a&-a&-ay(t+u)&-ay(t+u)&0\\
&&&&\\
\hline
&&&&\\
0&1&0&-ax(t+u)&0\\
&&&&\\
0&0&1-a(x+y)&-a(x+y)&-ay^2(t^2+tu+t^2)\\
&&&&\\
0&0&0&1-ax&-axy(t^2+tu+t^2)\\
\end{array}
} \right).
\]
Here our key result is as follows:
\begin{thm}\label{th:main1}
Let ${n}\geq1$ be a positive integer. Then we have
\begin{equation}
\frac{\det P_{n}}{\det P_{{n}-1}}
=(-1)^{{n}-1}ax^{{n}-1}[{n}]_{t,u}\prod_{i=1}^{{n}-1}(1-ax^{i}[{n}-i]_{x,y}),
\label{eq:Pn}
\end{equation}
and
\begin{equation}
\frac{\det P_{n}^{k}}{\det P_{n}}
=a\,x^{\frac{({k}-1)({k}-2)}{2}-\frac{{n}({n}-1)}{2}}\,
y^{\frac{({n}+1-{k})({n}+2-{k})}{2}}
\,[{n}+1]_{t,u} \left[{{{n}+1}\atop{{k}-1}}\right]_{x,y}
\label{eq:Pnk}
\end{equation}
for $1\leq {k}\leq {n}$,
\begin{equation}
\frac{\det P_{n}^{n+1}}{\det P_{n}} =a\,y\,[n+1]_{t,u}\,[n]_{x,y}
\end{equation}
and $\det P_{n}^{n+2}=0$.
\end{thm}

We need the
following:
\begin{lem}\label{lem:key}
For $0\leq m\leq {n}$,
\begin{align}
&\sum_{k=0}^{m} (-1)^{m-k}
x^{k\choose 2}y^{n-k\choose 2}
\left[{{n}\atop{k}}\right]_{x,y}
\prod_{i=0}^{k-1}\left\{1-ax^{i}[n-i]_{x,y}\right\}
\prod_{i=k}^{m-1}\left\{-ax^{i}[n-i]_{x,y}\right\}\nonumber\\
&\hspace{2cm}=x^{m\choose 2}y^{n-m\choose 2}
\left[{{n}\atop{m}}\right]_{x,y}
\prod_{i=1}^{m}\left\{1-ax^{i}[n-i]_{x,y}\right\}.\label{eq:matrixAC}
\end{align}
\end{lem}

\begin{demo}{Proof}
Note that
\[
[n]_{x,y}=y^{n-1}[n]_{x/y},\qquad {n\brack
k}_{x,y}=y^{kn-k^2}{n\brack k}_{x/y}.
\]
Since ${n-k\choose 2}+{k\choose 2}+kn-k^2={n\choose 2}$, setting
$c=ay^{n-1}$ and $q=x/y$, we can rewrite \eqref{eq:matrixAC} as
follows:
\begin{align}
\sum_{k=0}^{m} (-1)^{m-k} q^{k\choose 2} &{{n}\brack {k}}_{q}
\prod_{i=0}^{k-1}\left\{1-cq^{i}[n-i]_{q}\right\}
\prod_{i=k}^{m-1}\left\{-cq^{i}[n-i]_{q}\right\}\nonumber\\
=\,&q^{m\choose 2} {{n}\brack {m}}_{q}
\prod_{i=1}^{m}\left\{1-cq^{i}[n-i]_{q}\right\}.\label{eq:key}
\end{align}
Setting
\[
X=1+\frac{cq^n}{1-q},\quad Y=\frac{c}{1-q+cq^n},\quad
Z=\frac{cq^n}{1-q},
\]
then
$1-cq^i[n-i]_{q}=X(1-Yq^i)$ and $-cq^i[n-i]_{q}=Z(1-q^{i-n})$.
Hence, in \eqref{eq:key} making the following substitutions:
\begin{align*}
\prod_{i=0}^{k-1}\left\{1-cq^{i}[n-i]_{q}\right\} &=X^k(Y;q)_k,\\
\prod_{i=1}^{m}\left\{1-cq^{i}[n-i]_{q}\right\}&=X^m(Yq;q)_m,\\
\prod_{i=k}^{m-1}\left\{-cq^{i}[n-i]_{q}\right\}
&=(-1)^mZ^{m-k}q^{{m\choose
2}-mn}\frac{(q^{n-m+1};q)_m}{(q^{-n};q)_k},
\end{align*}
and writing the $q$-binomial coefficients as
\begin{align*}
{n\brack k}_{q}
=(-1)^kq^{kn-k(k-1)/2}\frac{(q^{-n};q)_k}{(q;q)_{k}},\quad {n\brack
m}_{q}=(-1)^mq^{mn-m(m-1)/2}\frac{(q^{-n};q)_m}{(q;q)_{m}},
\end{align*}
we see, after simplifying, that
identity~\eqref{eq:key} is equivalent to the special case
$Y=c/(1-q+cq^n)$ of the  identity:
$$
\sum_{k=0}^m\frac{(Y;q)_k}{(q;q)_k}Y^{m-k}=\frac{(Yq;q)_m}{(q;q)_m},
$$
which can be easily verified by induction.
\end{demo}

\begin{demo}{Proof of Theorem \ref{th:main1}}
We proceed by induction on ${n}$. When ${n}=1$, by a direct
computation we obtain $\det P_{1}=a$, $\det P_{1}^{1}=\det
P_{1}^{2}=a^2\,y\,[2]_{t,u}$ and $\det P_{1}^{3}=0$. This shows the
theorem is true when ${n}=1$. Let ${n}$ be an integer $\geq2$.
Assume the theorem is true for ${n}-1$.
\begin{enumerate}
\item[(i)]\label{roster:1}
We get
\begin{equation*}
\det P_{n}=\det\left(
\begin{array}{c|c}
P_{{n}-1}&\Pa_{{n}-1}\\
\noalign{\medskip} \hline\noalign{\medskip}
 X_{{n}-1}&\Pb_{{n}-1}
\end{array}\right)
=\det P_{{n}-1}\cdot
\det\left(\Pb_{{n}-1}-X_{{n}-1}P_{{n}-1}^{-1}\Pa_{{n}-1}\right)
\end{equation*}
and
\begin{equation*}
\det P_{n}^{k}=\det\left(
\begin{array}{c|c}
P_{{n}-1}&\Pa_{{n}-1}\\
\noalign{\medskip}\hline\noalign{\medskip} X_{{n}-1}&\Pb_{{n}-1}^k
\end{array}\right)
=\det P_{{n}-1}\cdot
\det\left(\Pb_{{n}-1}^{k}-X_{{n}-1}P_{{n}-1}^{-1}\Pa_{{n}-1}\right).
\end{equation*}
\item[(ii)]\label{roster:2}
By direct computation we can see that the $(i,j)$th entry of
$X_{{n}-1}P_{{n}-1}^{-1}\Pa_{{n}-1}$ ($1\leq i,j\leq n+1$) is equal
to
\[
\begin{cases}
\frac{\det P_{{n}-1}^{j}}{\det P_{{n}-1}}
&\text{ if $i=1$,}\\
0 &\text{ otherwise.}
\end{cases}
\]

By the induction hypothesis, the $(1,j)$th entry of
$X_{{n}-1}P_{{n}-1}^{-1}\Pa_{{n}-1}$ equals
\begin{equation}
\begin{cases}
a\,x^{{j-1\choose 2}-{n-1\choose 2}}\,y^{{n}+1-j\choose 2}\,[n]_{t,u}
\left[{{n}\atop{j-1}}\right]_{x,y}
&\text{ if $1\leq j\leq {n}-1$,}\\
a\,y\,[{n}-1]_{x,y}\,[n]_{t,u}
&\text{ if $j={n}$,}\\
0 &\text{ if $j={n}+1$.}
\end{cases}
\label{eq:entry}
\end{equation}
\item[(iii)]
Put
$W_{{n}-1}^{k}=\Pb_{{n}-1}^{k}-X_{{n}-1}P_{{n}-1}^{-1}\Pa_{{n}-1}$
and $W_{{n}-1}=\Pb_{{n}-1}-X_{{n}-1}P_{{n}-1}^{-1}\Pa_{{n}-1}$.
Then, by (i), we have $\frac{\det P_{n}^{k}}{\det P_{{n}-1}} =\det
W_{{n}-1}^{k} $ and $\frac{\det P_{n}}{\det P_{{n}-1}} =\det
W_{{n}-1}. $
By \thetag{\ref{eq:UM}} and \thetag{\ref{eq:entry}},
we can see that the $(1,j)$th
entry of $W_{{n}-1}^{k}$ is
\[
-a\,x^{{j-1\choose 2}-{{n}-1\choose 2}}\,y^{{{n}+1-j\choose 2}}[n]_{t,u}
\left[{{n}\atop{j-1}}\right]_{x,y}
\]
for $1\leq j\leq {n}$, and the $(1,{n}+1)$th entry is $0$ (the top
row does not depend on $k$). It is also easy to see that the
$(1,j)$th entry of $W_{{n}-1}$ is
\[
-a\,x^{{j-1\choose 2}-{{n}-1\choose 2}}\,y^{{{n}+1-j\choose 2}}[n]_{t,u}
\left[{{n}\atop{j-1}}\right]_{x,y}
\]
for $1\leq j\leq {n}+1$.
\item[(iv)]
We claim that
\[
\det
W_{{n}-1}=(-1)^{{n}-1}ax^{{n}-1}[n]_{t,u}\prod_{i=1}^{{n}-1}(1-ax^{i}[{n}-i]_{x,y}).
\]
In fact, the $(i,j)$th entry of $W_{{n}-1}$ is
\[
\begin{cases}
-a\,y^{\frac{({n}-j)({n}-j+1)}{2}}\,x^{\frac{(j-1)(j-2)}{2}-\frac{({n}-1)({n}-2)}{2}}[n]_{t,u}
\left[{{n}\atop{j-1}}\right]_{x,y}
&\text{ if $i=1$ and $1\leq j\leq {n}+1$,}\\
1-ax^{j-1}[{n}+1-j]_{x,y}
&\text{ if $i=j+1$ and $1\leq j\leq {n}$,}\\
-ax^{j-2}[{n}+2-j]_{x,y}
&\text{ if $i=j$ and $2\leq j\leq {n}+1$,}\\
0 &\text{ otherwise.}
\end{cases}
\]
Thus, if we expand $\det W_{{n}-1}$ along the top row, then we
obtain
\begin{align*}
\det W_{{n}-1}&=-ax^{-\frac{({n}-1)({n}-2)}2}[n]_{t,u}\\
&\quad \times \sum_{j=1}^{n+1}(-1)^{j+1}y^{\frac{(n-j)(n-j+1)}{2}}x^{\frac{(j-1)(j-2)}2}
\left[{{n}\atop{j-1}}\right]_{x,y}\det W_{n-1}(1;j).
\end{align*}
If we use
\[
\det W_{n-1}(1;j) =\prod_{\nu=0}^{j-2}(1-ax^{\nu}[n-\nu]_{x,y})
\prod_{\nu=j-1}^{n-1}(-ax^{\nu}[n-\nu]_{x,y}),
\]
then we obtain
\begin{align*}
\det W_{n-1} &=-ax^{-\frac{(n-1)(n-2)}2}[n]_{t,u}
\sum_{j=1}^{n+1}(-1)^{j+1}
y^{\frac{(n-j)(n-j+1)}{2}}x^{\frac{(j-1)(j-2)}2}\left[{{n}\atop{j-1}}\right]_{x,y}\\
&\hspace{2cm}\times \prod_{\nu=0}^{j-2}(1-ax^{\nu}[n-\nu]_{x,y})
\prod_{\nu=j-1}^{n-1}(-ax^{\nu}[n-\nu]_{x,y})\\
&=(-1)^{n-1}ax^{n-1}[n]_{t,u}\prod_{i=1}^{n-1}(1-ax^{i}[n-i]_{x,y})
\end{align*}
by \thetag{\ref{eq:matrixAC}}. Thus, by (i), we conclude that
\begin{equation}
\frac{\det P_{n}}{\det P_{n-1}} =\det W_{n-1}
=(-1)^{n-1}ax^{n-1}[n]_{t,u}\prod_{i=1}^{n-1}(1-ax^{i}[n-i]_{x,y}).
\label{eq:result1}
\end{equation}
\item[(v)]
We claim that
\[
\frac{\det P_{n}^{k}}{\det P_{n}} =\frac{\det W_{n-1}^{k}}{\det
W_{n-1}}=a\,x^{\frac{(k-1)(k-2)}{2}-\frac{n(n-1)}{2}}\,y^{\frac{(n+1-k)(n+2-k)}{2}}
\,[n+1]_{t,u} \left[{{n+1}\atop{k-1}}\right]_{x,y}
\]
for $1\leq k\leq n$. Because the rightmost column of $\Pb_{n-1}^{k}$
is the $k$th column of $U _{n}$, we have the $(i,n+1)$th entry of
$\Pb_{n-1}^{k}$ is
\[
\begin{cases}
-ay^n[n+1]_{t,u}
&\text{ if $i=2$,}\\
0 &\text{ otherwise,}
\end{cases}
\]
when $k=1$,
\[
\begin{cases}
-ay^{n+2-k}x^{k-2}[n+1]_{t,u}
&\text{ if $i=k$,}\\
-ay^{n+1-k}x^{k-1}[n+1]_{t,u}
&\text{ if $i=k+1$,}\\
0 &\text{ otherwise,}
\end{cases}
\]
when $2\leq k\leq n$,
\[
\begin{cases}
-ayx^{n-1}[n+1]_{t,u}
&\text{ if $i=n+1$,}\\
0 &\text{ otherwise,}
\end{cases}
\]
when $k=n+1$, and all zero when $k=n+2$. By the induction
hypothesis, the $(1,n+1)$th entry of $X_{n-1}P_{n-1}^{-1}\Pa_{n-1}$
is $ \frac{\det P_{n-1}^{n+1}}{\det P_{n-1}}=0. $
Thus the $(n+1)$th
column of $W_{n-1}^{k}=\Pb_{{n}-1}^{k}-X_{{n}-1}P_{{n}-1}^{-1}\Pa_{{n}-1}$ equals the $(n+1)$th column of
$\Pb_{n-1}^{k}$.
\begin{enumerate}
\item[(a)]
When $k=1$, we expand $\det W_{n-1}^{1}$ along the $(n+1)$th column,
then, by direct computation, we obtain
\[
\det W_{n-1}^{1} =(-1)^{n+3}(-ay^n[n+1]_{t,u})\det W_{n-1}(2;n+1)
\]
By expanding $\det W_{n-1}^{1}(2;n+1)$ along the top row we obtain
\[
\det
W_{n-1}(2;n+1)=\left(-a\,y^{\frac{n(n-1)}{2}}\,x^{-\frac{(n-1)(n-2)}{2}}[n]_{t,u}
\right)\prod_{\nu=1}^{n-1}(1-ax^{\nu}[n-\nu]_{x,y}).
\]
Thus we conclude that
\begin{align}
\det W_{n-1}^{1} &=(-1)^{n-1}a^2 y^{\frac{n(n+1)}{2}}
     x^{-\frac{(n-1)(n-2)}{2}}[n]_{t,u}\nonumber\\
&\qquad \qquad \times [n+1]_{t,u}\,\prod_{\nu=1}^{n-1}(1-ax^{\nu}[n-\nu]_{x,y}).
\end{align}
By \thetag{\ref{eq:result1}}, this implies
\[
\frac{\det W_{n-1}^{1}}{\det W_{n-1}} =ay^{\frac{n(n+1)}{2}}
x^{-\frac{n(n-1)}{2}}[n+1]_{t,u}
\]
which is the desired identity.
\item[(b)]
When $2\leq k\leq n$, we expand $\det W_{n-1}^{k}$ along the
$(n+1)$th column, then we obtain
\begin{align}
\det W_{n-1}^{k} &=(-1)^{k+n+1}(-ay^{n+2-k}x^{k-2}[n+1]_{t,u})\det
W_{n-1}(k;n+1)\nonumber\\
&+(-1)^{k+n+2}(-ay^{n+1-k}x^{k-1}[n+1]_{t,u})\det W_{n-1}(k+1;n+1).
\label{eq:expand}
\end{align}
By expanding along the top row, we obtain
\begin{align*}
\det W_{n-1}(k;n+1)
&=-ax^{-\frac{(n-1)(n-2)}{2}}[n]_{t,u}\,\sum_{j=1}^{k-1}(-1)^{j+1}
y^{\frac{(n-j)(n-j+1)}{2}}\,x^{\frac{(j-1)(j-2)}{2}}\\
&\times \left[{{n}\atop{j-1}}\right]_{x,y}\det W_{n-1}(1,k;j,n+1),\\
\det W_{n-1}(k+1;n+1)
&=-ax^{-\frac{(n-1)(n-2)}{2}}[n]_{t,u}\,\sum_{j=1}^{k}(-1)^{j+1}
y^{\frac{(n-j)(n-j+1)}{2}}
\,x^{\frac{(j-1)(j-2)}{2}}\\
&\times \left[{{n}\atop{j-1}}\right]_{x,y} \det
W_{n-1}(1,k+1;j,n+1),
\end{align*}
where $W_{n-1}(1,k;j,n+1)=W_{n-1}(k;n+1)(1;j)$ and
$W_{n-1}(1,k+1;j,n+1)=W_{n-1}(k+1;n+1)(1;j)$.
If we use
\begin{align*}
\det W_{n-1}(1,k;j,n+1)
&=\prod_{\nu=0}^{j-2}(1-ax^{\nu}[n-\nu]_{x,y})
\prod_{\nu=j-1}^{k-3}(-ax^{\nu}[n-\nu]_{x,y})\\
&\times\prod_{\nu=k-1}^{n-1}(1-ax^{\nu}[n-\nu]_{x,y}),\\
\det W_{n-1}(1,k+1;j,n+1)&=\prod_{\nu=0}^{j-2}(1-ax^{\nu}[n-\nu]_{x,y})
\prod_{\nu=j-1}^{k-2}(-ax^{\nu}[n-\nu]_{x,y})\\
&\times\prod_{\nu=k}^{n-1}(1-ax^{\nu}[n-\nu]_{x,y}),
\end{align*}
then, by \thetag{\ref{eq:matrixAC}}, we obtain
\begin{align*}
\det W_{n-1}(k;n+1) &=(-1)^{k-1}a
x^{\frac{(k-2)(k-3)}2-\frac{(n-1)(n-2)}2}
y^{\frac{(n-k+2)(n-k+1)}2} \\
&\times [n]_{t,u} \left[{{n}\atop{k-2}}\right]_{x,y}
\prod_{\nu=1}^{n-1}(1-ax^{\nu}[n-\nu]_{x,y}),
\\
\det W_{n-1}(k+1;n+1)
&=(-1)^{k}ax^{\frac{(k-1)(k-2)}2-\frac{(n-1)(n-2)}2}
y^{\frac{(n-k)(n-k+1)}2}\\
&\times [n]_{t,u} \left[{{n}\atop{k-1}}\right]_{x,y}
\prod_{\nu=1}^{n-1}(1-ax^{\nu}[n-\nu]_{x,y}).
\end{align*}
Thus, from \thetag{\ref{eq:expand}}, we conclude that
\begin{align*}
\det W_{n-1}^{k}
 =&(-1)^{n-1}a^2 x^{-\frac{(n-1)(n-2)}2+\frac{(k-1)(k-2)}2}
 y^{\frac{(n+1-k)(n+2-k)}2}
[n]_{t,u} [n+1]_{t,u}\\
&\times\left(y^{n-k+2}\left[{{n}\atop{k-2}}\right]_{x,y} +
x^{k-1}\left[{{n}\atop{k-1}}\right]_{x,y}\right)
\prod_{\nu=1}^{n-1}(1-ax^{\nu}[n-\nu]_{x,y})\\
& =(-1)^{n-1}a^2 x^{-\frac{(n-1)(n-2)}2+\frac{(k-1)(k-2)}2}
y^{\frac{(n+1-k)(n+2-k)}2} \\
&\quad \times [n]_{t,u} [n+1]_{t,u}
\left[{{n+1}\atop{k-1}}\right]_{x,y}
\prod_{\nu=1}^{n-1}(1-ax^{\nu}[n-\nu]_{x,y}).
\end{align*}
Using \thetag{\ref{eq:result1}}, we obtain
\[
\frac{\det W_{n-1}^{k}}{\det W_{n-1}}
=ax^{-\frac{n(n-1)}2+\frac{(k-1)(k-2)}2} y^{\frac{(n+1-k)(n+2-k)}2}
[n+1]_{t,u} \left[{{n+1}\atop{k-1}}\right]_{x,y},
\]
which is the desired identity.
\item[(c)]
When $k=n+1$, we also expand $\det W_{n-1}^{k}$ along the $(n+1)$th
column and repeat the same argument. It is not hard to obtain
\begin{align*}
\frac{\det W_{n-1}^{n+1}}{\det W_{n-1}} &=ay[n+1]_{t,u}[n]_{x,y}.
\end{align*}
The details are left to the reader.
\item[(d)]
When $k=n+2$,
$\det\Pb_{n-1}^{k}$ vanishes
since all the entries of the last column of $\det\Pb_{n-1}^{k}$ are zero.
\end{enumerate}
\end{enumerate}
This proves the theorem is true for $n$. By induction we conclude
that the theorem is true for all $n\geq1$. This completes the proof.
\end{demo}

Since $\det(P_1)=a$, Theorem \ref{thm:minor1}\;\eqref{eq:minor1}
follows easily from \eqref{eq:Pn}. \qed
\subsection{Proof of Theorem~\ref{thm:minor2}}

\bigbreak
Let $\seqpoly=\{\seqpoly_{n}\}_{n=1}^{\infty}$ be a sequence of
non-zero functions in finitely many variables $v_1,v_2,\dots$. We
use the convention that $\seqpoly_{n}!=\prod_{k=1}^{n}\seqpoly_{k}$
and
\[
\seqbinom{n}{k}=\begin{cases}
\frac{\seqpoly_{n}!}{\seqpoly_{k}!\seqpoly_{n-k}!},
&\text{ if $0\leq{k}\leq{n}$,}\\
0, &\text{ otherwise.}
\end{cases}
\]
We prove Theorem~\ref{thm:minor2}\;\eqref{eq:minor2} by considering the
following matrix $N_{n}(x,a)$, which generalize the matrix $N_n$
(set $x=1$ and $F_n=[n]_{t,u}$ to obtain $N_n$). Let $N_{n}({x},a)$
be the matrix defined inductively as follows:
\begin{align*}
&N_{0}({x},a)=\left( {x}\right)
\end{align*}
and
\begin{equation}
N_{n}({x},a)=\left(
\begin{array}{c|c}
N_{{n}-1}({x},a) &    \UB_{{n}-1}({x},a)\\
\noalign{\medskip} \hline \noalign{\medskip} O_{{n}+1,\N{{n}-1}} &
\LB_{{n}-1}({x},a)
\end{array}
\right) \label{eq:recB}
\end{equation}
where $\LB_{{n}-1}({x},a)$ is the $({n}+1)\times ({n}+1)$ matrix
defined by
\begin{equation}
\LB_{n-1}({x},a)=\left(\, {x}\delta_{ij}-a q^{i-1}[{n}+1-i]_{q}
(\delta_{ij}+\delta_{i+1,j}) \,\right)_{1\leq i,j\leq {n}+1}
\label{eq:LB}
\end{equation}
and $\UB_{{n}-1}({x},a)$ is the $\N{{n}-1}\times({n}+1)$ matrix
\[
\left(\begin{array}{c}
O_{\N{{n}-2},{n}+1}\\
\noalign{\medskip} \hline \noalign{\medskip} \UUB_{{n}-1}
\end{array}\right)
\]
with the ${n}\times({n}+1)$ matrix
\begin{equation}
\UUB_{{n}-1}= \left( -a\seqpoly_{n}\cdot(\delta_{ij}+\delta_{i+1,j})
\right)_{1\leq i\leq {n},\,1\leq j\leq{n}+1}.
\label{eq:UB}
\end{equation}
For instance, we get
\[
N_{2}({x},a)=\left( \scriptsize{
\begin{array}{cccccc}
{x}&-a\seqpoly_{1}&-a\seqpoly_{1}&0&0&0\\
&&&&&\\
0&{x}-a&-a&-a\seqpoly_{2}&-a\seqpoly_{2}&0\\
&&&&&\\
0&0&{x}&0&-a\seqpoly_{2}&-a\seqpoly_{2}\\
&&&&&\\
0&0&0&{x}-a(1+q)&-a(1+q)&0\\
&&&&&\\
0&0&0&0&{x}-aq&-aq\\
&&&&&\\
0&0&0&0&0&{x}\\
\end{array}
} \right).
\]

Let $\DB_{n}({x},a)$ denote the matrix obtained from $N_{n}({x},a)$
by deleting the $\N{n}$th row and the first column. Then the
following theorem is sufficient to prove our result. Here our
strategy is as follows. We regard $\det N_{n}({x},a)$ as a
polynomial in ${x}$ and find all linear factors. Finally we check
the leading coefficient in the both sides.
\begin{thm}
\label{th:conj} We have
\begin{equation}
\det \DB_{n}({x},a)
=(-1)^{\frac{{n}({n}-1)}2}a^{n}\,\seqpoly_{n}!\,{x}^{n}
\prod_{m=1}^{n-1}\prod_{k=1}^{n-m}\left({x}-aq^{k-1}[m]_{q}\right).
\label{eq:conj}
\end{equation}
\end{thm}

Then by setting $x=1$ and $F_n=[n]_{t,u}$ we obtain Theorem
\ref{thm:minor2}\;\eqref{eq:minor2}.\\

For instance, we have
\begin{align*}
\det \DB_{1}({x},a)&=\det\left(
\begin{array}{cc}
-a\seqpoly_{1}&-a\seqpoly_{1}\\
&\\
{x}-a&-a\\
\end{array}
\right)
={a}\,\seqpoly_{1}\,{x}\\
\intertext{and}
\det \DB_{2}({x},a)&=\det\left(
\begin{array}{ccccc}
-a\seqpoly_{1}&-a\seqpoly_{n}&0&0&0\\
&&&&\\
{x}-a&-a&-a\seqpoly_{2}&-a\seqpoly_{2}&0\\
&&&&\\
0&{x}&0&-a\seqpoly_{2}&-a\seqpoly_{2}\\
&&&&\\
0&0&{x}-a(1+q)&-a(1+q)&0\\
&&&&\\
0&0&0&{x}-aq&-aq\\
&&&&
\end{array}
\right)\\
&=-{a}^2\,\seqpoly_{1}\seqpoly_{2}\,{x}^2({x}-a).
\end{align*}
%
%

Fix positive integers $m$ and $k$. Define the row vectors
$\vecX_{n}^{{m},{k}}({x},{a})$ of degree $\N{n}$ as follows: For
$1\leq i\leq {n}+1$ and $1\leq j\leq i$, the
$\left(\frac{i(i-1)}2+j\right)$th entry of
$\vecX_{n}^{{m},{k}}({x},{a})$ is equal to
\begin{equation}
\vecX^{{m},{k}}_{i,j}= (-1)^{i+{m}+{k}}{a}
q^{-({m}+{k}-1)(i-{m}-{k})+\binom{j-{k}}{2}}
\frac{\seqpoly_{i-{m}-{k}}!}{[i-{m}-{k}]_{q}!}
\seqbinom{i-1}{{m}+{k}-1}\qbinom{m}{j-{k}}. \label{eq:eigen_vec}
\end{equation}
Here we use the convention that $\seqpoly_{n}!=[n]_{q}!=1$ if
$n\leq0$. For example, if ${n}=3$, ${m}={k}=1$, then
\begin{align*}
\vecX_{3}^{1,1}({x},{a})
=\left(
0,1,1,-\frac{\seqpoly_{2}}{q},-\frac{\seqpoly_{2}}{q},0,\frac{\seqpoly_{2}\seqpoly_{3}}{q^2[2]_{q}!},\frac{\seqpoly_{2}\seqpoly_{3}}{q^2[2]_{q}!},0,0
\right).
\end{align*}
\begin{lemma}
\label{lem:key2} Let ${n}$ be a positive integer. Let ${m}$ and
${k}$ be positive integers such that $1\leq{m}\leq{n}-1$ and
$1\leq{k}\leq{n}-{m}$. Then we have
\begin{equation}
\vecX_{n}^{{m},{k}}({x},{a})\, N_{n}({x},{a})
=({x}-{a}q^{{k}-1}[m]_{q})\,\vecX_{n}^{{m},{k}}({x},{a}).
\label{eq:key2}
\end{equation}
\end{lemma}
Before we proceed to the proof of the lemma, we see it in an
example. If ${n}=2$ and ${m}={k}=1$, then we have
\begin{align*}
&\left( 0,1,1,-\frac{\seqpoly_{2}}{q},-\frac{\seqpoly_{2}}{q},0
\right) \left( \scriptsize{
\begin{array}{cccccc}
{x}&-\seqpoly_{1}&-\seqpoly_{1}&0&0&0\\
&&&&&\\
0&{x}-a&-a&-a\seqpoly_{2}&-a\seqpoly_{2}&0\\
&&&&&\\
0&0&{x}&0&-a\seqpoly_{2}&-a\seqpoly_{2}\\
&&&&&\\
0&0&0&{x}-a(1+q)&-a(1+q)&0\\
&&&&&\\
0&0&0&0&{x}-aq&-aq\\
&&&&&\\
0&0&0&0&0&{x}\\
\end{array}
}
\right)\\
&= \left(
0,{x}-a,{x}-a,-\frac{\seqpoly_{2}}{q}({x}-{a}),-\frac{\seqpoly_{2}}{q}({x}-{a}),0
\right).
\end{align*}
\begin{demo}{Proof of Lemma~\ref{lem:key2}}
We proceed by induction on ${n}$. When ${n}=0$ or ${n}=1$, our claim
is easy to check by direct computation. Assume
\thetag{\ref{eq:key2}} is true upto ${n}-1$. Then the first
$\N{{n}-1}$ entries of $ \vecX_{n}^{{m},{k}}({x},{a}) N_{n}({x},{a})
$ agree with those of
$({x}-q^{{k}-1}[{m}]_{q})\vecX_{n}^{{m},{k}}({x},{a})$ by the
induction hypothesis. So we have to check the last ${n}+1$ entries.
In fact we verify the following three cases.
\begin{enumerate}
\item[(i)]
If $i={n}+1$ and $j=1$, then the
$\left(\frac{{n}({n}+1)}2+1\right)$th entry of
$\vecX_{n}^{{m},{k}}({x},{a})N_{n}({x},{a})$ is equal to
\begin{align*}
&(-a\,\seqpoly_{n})\vecX^{{m},{k}}_{n,1}
+({x}-{a}[{n}]_{q})\vecX^{{m},{k}}_{n+1,1}
.
\end{align*}
Note that the coefficient $\qbinom{m}{1-{k}}$ becomes zero unless
${k}=1$. Thus, by direct computation, one can easily check that this
sum equals
\begin{align*}
(-1)^{{n}+{m}+{k}+1}&{a}q^{-({m}+{k}-1)({n}-{m}-{k}+1)+\binom{1-{k}}2}\\
&\times
\frac{\seqpoly_{{n}+1-{m}-{k}}!}{[{n}+1-{m}-{k}]_{q}!}
\seqbinom{n}{{m}+{k}-1}\qbinom{m}{1-{k}}({x}-{a}[{m}+{k}-1]_{q}).
\end{align*}
\item[(ii)]
If $i={n}+1$ and $2\leq j\leq{n}$, then the
$\left(\frac{{n}({n}+1)}2+j\right)$th entry of
$\vecX_{n}^{{m},{k}}({x},a)N_{n}({x},{a})$ is equal to
\begin{align*}
&(-a\seqpoly_{n})\vecX^{{m},{k}}_{n,j-1}
+(-a\seqpoly_{n})\vecX^{{m},{k}}_{n,j}\\
&\qquad\qquad+\left(-{a}q^{j-2}[{n}-j+2]_{q}\right)\vecX^{{m},{k}}_{n+1,j-1}
+\left({x}-{a}q^{j-1}[{n}-j+1]_{q}\right)\vecX^{{m},{k}}_{n+1,j}
.
\end{align*}
By direct computation, one can easily check this equals
\begin{align*}
(-1)^{{n}+{m}+{k}+1}&{a}q^{-({m}+{k}-1)({n}-{m}-{k}+1)+\binom{j-{k}}2}\\
&\times
\frac{\seqpoly_{{n}+1-{m}-{k}}!}{[{n}+1-{m}-{k}]_{q}!}
\seqbinom{n}{{m}+{k}-1}\qbinom{m}{j-{k}}\left({x}-{a}q^{{k}-1}[{m}]_{q}\right).
\end{align*}
\item[(iii)]
If $i=j={n}+1$, then the $\frac{({n}+1)({n}+2)}2$th entry of
$\vecX_{n}^{{m},{k}}({x},{a})N_{n}({x},{a})$ is equal to
\begin{align*}
&(-a\seqpoly_{n})\vecX^{{m},{k}}_{n,n}
+\left(-{a}q^{n-1}\right)\vecX^{{m},{k}}_{n+1,n}
+{x}\vecX^{{m},{k}}_{n+1,n+1}
.
\end{align*}
One can easily check this is always equal to zero.
\end{enumerate}
Thus this completes the proof of our lemma.
\end{demo}
\begin{corollary}
\label{cor:factors} Let $n$ be a positive integer. Then there exists
a polynomial $\varphi({x})$ such that
\begin{equation*}
\det \DB_{n}({x},a)
=\varphi({x})\prod_{m=1}^{{n}-1}\prod_{k=1}^{n-m}\left({x}-aq^{k-1}[m]_{q}\right).
\end{equation*}
\end{corollary}
\begin{demo}{Proof}
Let $\vecU_{n}^{{m},{k}}({x},{a})$ (resp.
$\vecV_{n}^{{m},{k}}({x},{a})$) denote the vector of degree
$\N{n}-1$ obtained from $\vecX_{n}^{{m},{k}}({x},{a})$ by deleting
the last (resp. first) entry. Then, by \thetag{\ref{eq:key2}}, we
obtain
\[
\vecU_{n}^{{m},{k}}({x},{a})\DB_{n}({x},{a})=\left({x}-aq^{{k}-1}[{m}]_q\right)\,\vecV_{n}^{{m},{k}}({x},{a}).
\]
By substituting ${x}=q^{{k}-1}[{m}]$ into this identity we obtain
\[
\vecU_{n}^{{m},{k}}(aq^{{k}-1}[{m}]_q,{a})\DB_{n}(aq^{{k}-1}[{m}]_q,{a})=\vec0.
\]
Since $\vecU_{n}^{{m},{k}}({x},{a})$  is a non-zero vector when
$1\leq m\leq n-1$ and $1\leq k\leq n-m$, $\DB_{n}(aq^{p-1}[m]_q,a)$ is
singular, which means $\det \DB_{n}(aq^{{k}-1}[{m}]_q,{a})=0$. Thus we
conclude that $\det \DB_{n}({x},{a})$ is divisible by
${x}-aq^{{k}-1}[{m}]_q$, which immediately implies our corollary.
\end{demo}
\begin{prop}
Let $n$ be a positive integer. Then there exists a polynomial
$\psi({x})$ such that
\begin{equation*}
\det \DB_{n}({x},a)
=\psi({x})\,{x}^{n}\prod_{{m}=1}^{{n}-1}
\prod_{{k}=1}^{{n}-{m}}\left({x}-aq^{{k}-1}[m]_q\right).
\end{equation*}
\end{prop}
\begin{demo}{Proof}
By Corollary~\ref{cor:factors}, we only need to show that $\det
\DB_{n}({x},a)$ is divisible by ${x}^{n}$. We show this by the
following column transformations on $\DB_{n}({x},a)$. First note
that $\DB_{n}({x},a)$ has $\N{n}-1$ columns. For each
$i=1,\dots,{n}$, let $\Col{i}$ denote the set of columns
$j=\frac{i(i+1)}2,\frac{i(i+1)}2+1,\dots,\frac{(i+1)(i+2)}2-1$. We
perform the following column transformations in each block
$\Col{i}$. For each $i$, we subtract the $\frac{i(i+1)}2$th column
from the $\left(\frac{i(i+1)}2+1\right)$th column, then subtract the
$\left(\frac{i(i+1)}2+1\right)$th column from the
$\left(\frac{i(i+1)}2+2\right)$th column, and so on, until we
subtract the $\left(\frac{(i+1)(i+2)}2-2\right)$th column from the
$\left(\frac{(i+1)(i+2)}2-1\right)$th column. Then each entry of the
$\left(\frac{(i+1)(i+2)}2-1\right)$th column becomes $\pm {x}$. For
example, if we perform this operation to $\DB_{3}({x},{a})$ which
looks like
\[
\left(
\begin{array}{cc|ccc|cccc}
 -{a}\seqpoly_{1}  & -{a}\seqpoly_{1}&           0&          0&          0&           0&              0&      0&     0\\
{x}-{a}& -{a}& -a\seqpoly_{2}&-a\seqpoly_{2}&          0&           0&              0&      0&     0\\
      0&  {x}&           0&-a\seqpoly_{2}&-a\seqpoly_{2}&           0&              0&      0&     0\\
      0&    0&{x}-a[2]_{q}&  -a[2]_{q}&          0& -a\seqpoly_{3}&  -{a}\seqpoly_{3}&      0&     0\\
      0&    0&           0&  {x}-{a}q &        -aq&           0&  -{a}\seqpoly_{3}&-a\seqpoly_{3}&     0\\
      0&    0&           0&          0&        {x}&           0&              0&-a\seqpoly_{3}& -a\seqpoly_{3}\\
      0&    0&           0&          0&          0&{x}-a[3]_{q}&    -{a}[3]_{q}&      0&     0\\
      0&    0&           0&          0&          0&           0&{x}-{a}q[2]_{q}& -aq[2]_{q}&     0\\
      0&    0&           0&          0&          0&           0&              0& {x}-aq^2& -aq^2\\
\end{array}\right),
\]
then we obtain
\[
\left(
\begin{array}{cc|ccc|cccc}
 -a\seqpoly_{1}&  0&     0&     0&     0&     0&      0&      0&     0\\
{x}-a& -{x}& -a\seqpoly_{2}&     0&     0&     0&      0&      0&     0\\
  0&  {x}&     0& -a\seqpoly_{2}&     0&     0&      0&      0&     0\\
  0&  0&{x}-a[2]_{q}&    -{x}&     {x}& -a\seqpoly_{3}&      0&      0&     0\\
  0&  0&     0& {x}-aq &    -{x}&     0&  -a\seqpoly_{3}&      0&     0\\
  0&  0&     0&     0&     {x}&     0&      0&  -a\seqpoly_{3}&     0\\
  0&  0&     0&     0&     0&{x}-a[3]_{q}&     -{x}&      {x}&    -{x}\\
  0&  0&     0&     0&     0&     0&{x}-aq[2]_{q}&     -{x}&     {x}\\
  0&  0&     0&     0&     0&     0&      0& {x}-aq^2&    -{x}\\
\end{array}\right).
\]
This is always true. One can easily check the last column of each
block becomes $\pm {x}$ after these elementary transformations using
the definition of $\LB_{n-1}({x},a)$ and $\UB_{n-1}({x},a)$ in
\thetag{\ref{eq:LB}} and \thetag{\ref{eq:UB}}. Thus, by taking the
determinant of $\DB_{n}({x},{a})$, we can factor out ${x}$ from each
block $\Col{i}$, ${i}=1,2,\dots,{n}$, and we conclude that
$\det\DB_{n}({x},{a})$ is divisible by ${x}^{n}$.
\end{demo}
Now we are in position to complete the proof of
Theorem~\ref{th:conj}.
\begin{demo}{Proof of Theorem~\ref{th:conj}}
To complete the proof of Theorem~\ref{th:conj}, we need to show that
the degree of $\det\DB_{n}({x},{a})$
 is $\frac{{n}({n}+1)}{2}$ as a polynomial in ${x}$,
and the leading coefficient of $\det\DB_{n}({x},{a})$ is equal to
$(-1)^{\frac{{n}({n}-1)}2}{a}^{n}\seqpoly_{n}!$. Let $\K{n}=\N{n}-1$
which is the degree of the matrix $\DB_{n}({x},a)$. Let
$\entryB_{ij}$ denote the $(i,j)$th entry of $\DB_{n}({x},{a})$. By
the definition of determinants we have
\[
\det \DB_{n}({x},{a}) =\sum_{\pi\in S_{\K{n}}} \sgn\pi\,
\entryB_{\pi(1)1}\entryB_{\pi(2)2}\cdots\entryB_{\pi(\K{n})\K{n}}.
\]
We use the two-line notation
\[
\pi=\left(\begin{array}{cccc}
1&2&\hdots&\K{n}\\
\pi(1)&\pi(2)&\hdots&\pi(\K{n})
\end{array}\right)
\]
to express a permutation $\pi$ of letters $[\K{n}]$. For each $j$,
if $\pi(j)=j+1$, then the entry $\entryB_{\pi(j)j}$ is of degree $1$
as a polynomial in $x$, and otherwise it is a constant. Thus $\det
\DB_{n}({x},a)$ is apparently of at most
$\K{n}-1=\frac{(n+1)(n+2)}{2}-2$ degree as a polynomial in ${x}$.
For example $\DB_{3}({x},a)$ looks as follows.
\[
\left(
\begin{array}{cc|ccc|cccc}
 -a\seqpoly_{1}& \pmb{-a\seqpoly_{1}}&     0&     0&     0&     0&      0&      0&     0\\\hline
 \pmb{{x}-a}& -a& -a\seqpoly_{2}& -a\seqpoly_{2}&     0&     0&      0&      0&     0\\
  0&  {x}&     0& -a\seqpoly_{2}&\pmb{-a\seqpoly_{2}}&     0&      0&      0&     0\\\hline
  0&  0&\pmb{{x}-a[2]_{q}}& -a[2]_{q}&     0& -a\seqpoly_{3}&  -a\seqpoly_{3}&      0&     0\\
  0&  0&     0&\pmb{{x}-aq}&   -aq&     0&  -a\seqpoly_{3}&  -a\seqpoly_{3}&     0\\
  0&  0&     0&     0&     {x}&     0&      0&  -a\seqpoly_{3}& \pmb{-a\seqpoly_{3}}\\\hline
  0&  0&     0&     0&     0&\pmb{{x}-a[3]_{q}}&  -a[3]_{q}&      0&     0\\
  0&  0&     0&     0&     0&     0&\pmb{{x}-aq[2]_{q}}& -aq[2]_{q}&     0\\
  0&  0&     0&     0&     0&     0&      0&\pmb{{x}-aq^2}& -aq^2\\
\end{array}\right),
\]
Our first claim is that $\det \DB_{n}({x},a)$ is a polynomial of
degree $\frac{n(n+1)}{2}$. Let $\Col{i}$, $i=1,2,\dots,{n}$, be as
in the previous proof. Note that $\Col{i}$ includes $i+1$ columns.
We claim that $\pi(j)=j+1$ can happen at most $i$ column indices $j$
in each block $\Col{i}$. Otherwise
$\entryB_{\pi(1)1}\entryB_{\pi(2)2}\cdots\entryB_{\pi(\K{n})\K{n}}$
vanishes. In fact, assume that $\pi(j)=j+1$ for all $j$ in a certain
block $\Col{i}$. Then this must be the case for the block
$\Col{i+1}$. There is no other choice if we assume
$\entryB_{\pi(1)1}\entryB_{\pi(2)2}\cdots\entryB_{\pi(\K{n})\K{n}}$
is nonzero. And this must be also the case for the block
$\Col{i+2}$, and so on. Finally we have to take  $\pi(j)=j+1$ for
all $j$ in the block $\Col{n}$, but this is impossible. Thus we
reach a contradiction. We conclude that the degree of
$\det\DB_{n}({x},{a})$ is at most $\frac{n(n+1)}2$. In fact there is
a permutation which realize this degree, i.e.
\[
\pi=\left(\begin{array}{cc|ccc|cc|ccccc}
1&2&3&4&5&\hdots&\hdots&\K{n-1}+1&\K{n-1}+2&\hdots&\K{n}-1&\K{n}\\
2&1&4&5&3&\hdots&\hdots&\K{n-1}+2&\K{n-1}+3&\hdots&\K{n}&\K{n-1}+1
\end{array}\right).
\]
It is easy to see that this $\pi$ is the only permutation with which
$\entryB_{\pi(1)1}\entryB_{\pi(2)2}\cdots\entryB_{\pi(\K{n})\K{n}}$
does not vanish and of degree $\frac{n(n+1)}2$. Thus we conclude
that the leading coefficient of $\det \DB_{n}({x},a)$ equals
\[
\sgn\pi\cdot (-1)^{n}\, a^{n}\, \seqpoly_{n}!.
\]
This immediately implies the resulting identity
\thetag{\ref{eq:conj}}.
\end{demo}

\begin{rmk}
One may notice that $M_{n}$ in \thetag{\ref{eq:recM}} and $N_{n}$ in
\thetag{\ref{eq:recB}} are in a similar form, but our methods to
evaluate them are far from parallel. It seems that the first method
does not work with the matrix $N_{n}$ since we can't guess the
entries of $CA^{-1}N$ as we did in \thetag{\ref{eq:Pnk}}. Meanwhile,
the second method does not work with the matrix $M_{n}$ at this
point since even if we generalize $M_{n}$ to $M_{n}({x},a)$, we
don't know the general form of the eigenvectors of $M_{n}({x},a)$.
The reader can find the general guidance about matrix evaluation in
\cite{Kr}. We may say that the second proof follows this general
philosophy.
\end{rmk}

\end{document}